\documentclass[draftcls,onecolumn]{IEEEtran}
\usepackage{xr-hyper, cite}
\usepackage{comment}
\usepackage{amsmath,amssymb,amsfonts}
\usepackage{multirow}
\usepackage{mathtools}
\usepackage{graphicx}
\usepackage{textcomp}
\usepackage{xcolor}
\usepackage{hyperref}
\usepackage{amsbsy}
\usepackage{amsthm, enumitem}
\usepackage{booktabs,tabularx}
\newcolumntype{Y}{>{\raggedright\arraybackslash}X}

\usepackage{algorithm,algpseudocode}
\algnewcommand{\Inputs}[1]{%
	\State \textbf{Inputs:}
	\Statex \hspace*{\algorithmicindent}\parbox[t]{.8\linewidth}{\raggedright #1}
}
\algnewcommand{\Initialize}[1]{%
	\State \textbf{Initialize:}
	\Statex \hspace*{\algorithmicindent}\parbox[t]{.8\linewidth}{\raggedright #1}
}
\makeatletter
\AtEndDocument{%
	\begingroup
	\edef\@currentlabel{\number\value{equation}}\label{ctr:last-equation}%
	\edef\@currentlabel{\number\value{lemma}}\label{ctr:last-lemma}%
    \edef\@currentlabel{\number\value{section}}\label{ctr:last-app}%
	\endgroup
}
\makeatother
\setlist[itemize]{leftmargin=*}
\usepackage{subcaption}
\usepackage{bm}

\providecommand{\abs}[1]{\lvert#1\rvert}
\newcommand{\br}[1]{\left({#1}\right)}

\newcommand{\bsq}[1]{\left[{#1}\right]}
\renewcommand{\b}[1]{\ensuremath{\mathbf{#1}}} 
\newcommand{\E}[1]{\ensuremath{\mathbb{E}\left[#1\right]}}  
\newcommand{\En}[1]{\ensuremath{\mathbb{E}#1}}  
\newcommand{\Ex}[1]{\ensuremath{\EE_{\xi}[#1]}}  
\newcommand{\Et}[1]{\ensuremath{\mathbb{E}_{i_t}[#1]}}  
\newcommand{\Ej}[1]{\ensuremath{\mathbb{E}_{j_t}[#1]}}  
\newcommand{\norm}[1]{\ensuremath{\left\|#1\right\|}} 
\newcommand{\eqtext}[1]{\ensuremath{\stackrel{\text{#1}}{=}}} 
\newcommand{\leqtext}[1]{\ensuremath{\stackrel{\text{#1}}{\leq}}} 
\newcommand{\geqtext}[1]{\ensuremath{\stackrel{\text{#1}}{\geq}}} 
\DeclareMathOperator*{\argmin}{arg\,min} 
\providecommand{\ip}[2]{\langle #1, #2 \rangle} 
\newcommand{\col}[1]{\textcolor{black}{#1}} 
\definecolor{bottlegreen}{RGB}{0, 106, 78} 

\renewcommand{\O}[1]{{\mathcal{O}\left(#1\right)}}
\newcommand{\Ot}[1]{{\tilde{\mathcal{O}}\left(#1\right)}}
\newcommand{\prox}[2]{{\text{prox}_{#1}\!\left(#2\right)}}

\def \x {{\b{x}}}
\def \y {{\b{y}}}
\def \z {{\b{z}}}

\def \a {{\b{a}}}

\def \v {{\b{v}}}
\def \w {{\b{w}}}
\def \u {{\b{u}}}

\def \EE {{\mathbb{E}}}

\def \lam {{\boldsymbol{\lambda}}}
\def \mub {{\boldsymbol{\mu}}}
\def \T {{\mathsf{T}}}

\def \cX {{\mathcal{X}}}

\def \cK {{\mathcal{K}}}

\def \cN {{\mathcal{N}}}
\def \cS {{\mathcal S}}

\def \sX {{\mathsf{X}}}
\def \sY {{\mathsf{Y}}}

\def \eps {{\epsilon}}

\def \Rn {{\mathbb{R}}}

\def \tx {{\tilde{\x}}}
\def \bx {{\bar{\x}}}

\def \bt {{\tilde{B}}}

\def \tg {{\tilde{g}}}
\def \tphi {{\tilde{\phi}}}

\newtheorem{assumption}{}

\theoremstyle{remark}

\newtheorem{theorem}{Theorem}
\newtheorem{lemma}{Lemma}

\definecolor{bottlegreen}{RGB}{0, 106, 78}

\begin{document}
	
	\title{Hinge-Proximal Stochastic Gradient Methods for Convex Optimization with Functional Constraints\author{Vaibhav Rajoriya,  Prateek Priyaranjan Pradhan, Ketan Rajawat \thanks{Authors are with the Dept. of Electrical Engineering, Indian Institute of Technology Kanpur, Kanpur, Uttar Pradesh 208016, India. }}}
	\maketitle
	
	\begin{abstract}
	This paper considers stochastic convex optimization problems with smooth functional constraints arising in constrained estimation and robust signal recovery. We operate in the high-dimensional and highly-constrained setting, where oracle access is restricted to one or a few objective and constraint gradients per-iteration, as in streaming or online estimation. Existing approaches to solve such problems are based on either the stochastic primal-dual or stochastic subgradient methods, and require globally Lipschitz continuous constraint functions. In this work, we develop a hinge-proximal framework that utilizes an exact penalty reformulation to yield updates involving only one linearized constraint (and hence accessing one constraint gradient) per-iteration. The updates also admit a novel hinge-proximal three-point inequality relying on smoothness rather than global Lipschitz continuity of the constraint functions. The framework leads to three algorithms: a baseline hinge-proximal SGD (HPS), a variance-reduced HPS version for finite-sum settings, and a nested HPS version whose performance depends on a geometric regularity constant of the constraint region rather than explicitly on the number of constraints, while achieving near-SGD sample complexity. The superior empirical performance of the proposed algorithms is demonstrated on a robust regression problem with noisy features, representative of errors-in-variables estimation.
	\end{abstract}
	
\section{Introduction}\label{sec:intro}
In this paper, we focus on solving the following constrained optimization problem
\begin{align}\tag{$\mathcal{P}$}\label{mainProb}
    \begin{aligned}
        \mathbf{x}_\star = \argmin_{\mathbf{x} \in \mathbb{R}^d} ~ &f(\mathbf{x}) + h(\mathbf{x}), \\
        \text{s. t. } \hspace{5mm} &g_j(\mathbf{x}) \leq 0 \hspace{1cm} 1\leq j\leq m
    \end{aligned}
\end{align}
where the objective function $f:\Rn^d \rightarrow \Rn$ is strongly convex with $f(\mathbf{x}) := \Et{f_{i_t}(\mathbf{x})}$ and $\mathbb{E}_{i_t}[\cdot]$ denotes the expectation with respect to the random index $i_t$. We also consider the finite-sum case where $f(\x) = \frac{1}{n}\sum_{i=1}^n f_i(\x)$ and $i_t$ sampled uniformly from $\{1, \ldots, n\}$. The functions $f_i: \mathbb{R}^d \rightarrow \mathbb{R}$ and $g_j: \mathbb{R}^d \rightarrow \mathbb{R}$ are proper, closed, and convex, and smooth. The regularization function $h: \mathbb{R}^d \rightarrow \mathbb{R}$ is convex but potentially non-smooth and may include an indicator function corresponding to a set-inclusive constraint of the form $\x \in \cK$ for closed convex set $\cK$. Stochastic and finite-sum optimization problems with non-linear inequality constraints arise in signal processing and related areas, such as constrained estimation and regression~\cite{zhang2015joint, liu2019stochastic}, fair classification~\cite{fair_icassp, akhtar2021conservative}, trajectory optimization~\cite{malyuta2022convex}, constrained federated learning~\cite{zhang2025federated}, constrained reinforcement learning~\cite{tian2022successive}, and optimal power flow~\cite{owerko2025learning}.

In many of these applications, the target regimes are simultaneously  high-dimensional (large $d$) and highly constrained (large $m$). Hence, direct access to $f$ or the full collection $\{g_j\}_{j=1}^m$ is typically out of reach, rendering projection-based methods impractical and motivating sample-based and oracle-efficient approaches. Of particular interest are algorithms designed to minimize the number of calls to the stochastic first-order oracle (SFO), which for a given $\x$, returns a stochastic gradient from a randomly selected measurement or sample $\nabla f_{i_t}(\x)$ together with one randomly sampled functional constraint $g_{j_t}(\x)$ and $\nabla g_{j_t}(\x)$. The oracle model is natural in streaming and online settings, where each update processes only a single data record and a single constraint function.

Existing first-order approaches for solving \eqref{mainProb} are based on two main paradigms: primal-dual algorithms \cite{xu2020primal, singh2025stochastic} and primal algorithms inspired from Polyak's subgradient iteration~\cite{nedic2019random, necoara2022stochastic, singh2024stochastic, singh2024stochastic_movingball, singh2024unified, singh2024mini}. Of these, \cite{xu2020primal, necoara2022stochastic} consider the stochastic version of \eqref{mainProb}, make only one SFO call per-iteration, and achieve SFO complexity comparable to unconstrained stochastic gradient descent (SGD). However, both of these methods require the constraints to be globally Lipschitz, i.e., have uniformly bounded subgradients, an assumption violated in many settings where smooth squared-error, quadratic, log-sum-exp, and barrier-type constraints are more natural. The finite-sum version of \eqref{mainProb} can be solved using \cite{singh2024mini, singh2024unified} or related deterministic algorithms~\cite{nedic2019random, singh2025stochastic, singh2024stochastic, singh2024stochastic_movingball}, but all of these achieve the suboptimal $\O{1/\epsilon}$ SFO complexity while also requiring globally Lipschitz constraints. Additionally, all the existing primal-domain algorithms need a linear regularity assumption on the feasible set. 

In this work, we put forth a new hinge-proximal framework for stochastic optimization with smooth functional constraints. The proposed framework utilizes an exact hinge-penalty reformulation to develop a class of primal-domain SGD algorithms for solving \eqref{mainProb} while using only one randomly sampled constraint gradient per-iteration. We instantiate this framework in three algorithms: the baseline hinge-proximal SGD (HPS), a variance-reduced HPS (VR-HPS) for finite-sum problems, and a nested HPS (N-HPS) that leverages linear regularity to remove the explicit dependence on the number of constraints $m$. A key technical ingredient of the framework is a hinge-proximal analogue of the classical three-point property for the proximal operator, which we use to analyze all three variants. The three algorithms have the following distinguishing features.
\begin{itemize}
\item The HPS algorithm makes one SFO call per-iteration and achieves an SFO complexity of $\O{m^2/\epsilon}$ without requiring any regularity condition on the constraints. While this dependence on $m$ is suboptimal, HPS serves as a conceptual and technical baseline for developing the other variants.
\item Building upon the hinge-proximal framework, we next develop the VR-HPS algorithm, which exploits the finite-sum structure of both the objective and the penalty, achieving an SFO complexity of \col{$\O{n+\tfrac{m+\sqrt{n}}{\sqrt{\eps}}}$} without linear regularity. Variance reduction in the non-smooth penalty term is obtained by adapting the stochastic decoupling idea of \cite{mishchenko2019stochastic} to hinge penalties and linearized constraints.
\item Finally, seeking to remove the explicit dependence on $m$ in the rate expressions, we develop the N-HPS algorithm, which solves the stochastic version of \eqref{mainProb} under the linear regularity condition. N-HPS uses an inner loop to (approximately) project onto a single randomly selected constraint via an exact-penalty subproblem, and achieves an SFO complexity of $\Ot{\tfrac{1}{\eps}}$ with no explicit dependence on $m$ in the leading term.
\end{itemize}
In summary, for solving \eqref{mainProb} with smooth constraint functions and without global Lipschitz (bounded-subgradient) assumptions, HPS and N-HPS are the first stochastic algorithms that achieve $\Ot{\tfrac{1}{\eps}}$ SFO complexity (see Table~\ref{tab:sc-constraints1}), while VR-HPS is the first algorithm that achieves $\O{\tfrac{1}{\sqrt{\eps}}}$ complexity for the finite-sum case with access to only one constraint gradient per-iteration (see Table~\ref{tab:sc-constraints2}).

\subsection{Robust Regression}
To motivate the usefulness and demonstrate the efficacy of the proposed algorithms, we evaluate its numerical performance on the problem of robust regression under noisy features. Recall that the general regression problem with a linear measurement model entails estimating $\x$ using 
\begin{equation}\label{ols}
    \min_{\x} \frac{1}{n}\sum_{i=1}^n \ell(\a_i^\top \x, b_i).   
\end{equation}
where $\ell:\Rn\times \Rn \rightarrow \Rn_+$ is a convex loss function. However, in many real world domains, such as in sensor networks, wireless communication, and biomedical signal analysis, the features are often corrupted by complex and poorly understood acquisition noise. Estimation of $\x$ from noisy feature-vector pairs $(\a_i,b_i)_{i=1}^{n}$ has been widely studied under the umbrella of errors-in-variables and measurement-error regression~\cite{carroll2006measurement}. 

Some of the optimization-based approaches for handling noisy features include worst-case robust optimization~\cite{bertsimas2011theory}, distributionally robust optimization~\cite{wang2021distributionally}, and total least-squares~\cite{markovsky2007overview}. These methods rely on explicit parametric descriptions of the uncertainty and can be overly conservative if the true corruption mechanism is complex or only available only through simulation. On the other hand, Bayesian approaches model measurement noise probabilistically and can flexibly specify parametric noise models, but they also require an explicit noise model while relying on computationally intensive interfaces like Markov Chain Monte Carlo. 

Here we consider an alternative framework that retains the loss function in the objective but introduces constraints that keep the per-sample loss small for sampled corrupted features: 
\begin{equation}
\begin{aligned}
    \min_{\x} \quad & \frac{1}{n}\sum_{i=1}^n \ell(\a_i^\top \x, b_i)  \\
    \text{s. t.} \quad & \ell(P_j(\a_i)^\top \x, b_i) \leq \varepsilon, \quad \forall i,j
\end{aligned}
\label{eqn:robust_regression}
\end{equation}
where $P_j(\a_i)$ denotes the perturbed realization of the true feature vector $\a_i$ and $\varepsilon$ is a user-defined tolerance for prediction error under perturbations. The constraints in \eqref{eqn:robust_regression} ensure that even under feature corruption, the prediction loss remains within an acceptable threshold of  $\varepsilon$ for the considered perturbation realizations. This formulation allows us to incorporate Monte-Carlo models of feature noise, arising from sensor simulators, calibration uncertainty, and adversarial perturbation generators, without an explicit parametric noise model. However, this formulation is less common in the literature because it introduces a large number of constraints. This motivates the proposed stochastic methods that access only a single data record and a single sampled corruption constraint per iteration. In the remainder of the paper, we develop and analyze such methods, and demonstrate their effectiveness on robust regression tasks.

\subsection{Related Works}
In this section, we discuss \cite{nedic2019random, necoara2022stochastic, singh2024unified, singh2024mini, xu2020primal, singh2025stochastic, singh2024stochastic, singh2024stochastic_movingball} and other related works in greater detail. For the sake of brevity, we only include works requiring first-order objective and constraint information, and exclude works involving projections onto the feasible set~\cite{yang2017richer}. 

Primal-dual algorithms for solving \eqref{mainProb} maintain and update dual variables and hence require $\O{m}$ storage per-iteration \cite{xu2020primal, singh2025stochastic}. Generally, these methods struggle to control the norm of the dual iterates, which is either assumed to be bounded \cite{singh2025stochastic} or carefully bounded through the choice of parameters \cite{xu2020primal}. For the stochastic version of \eqref{mainProb}, the SFO complexity of \cite{xu2020primal} is $\Ot{\frac{1}{\eps}}$ which is almost at par with that of SGD, up to a logarithmic factor. The proposed N-HPS also achieves the same complexity under a linear regularity assumption, but without requiring global Lipschitz bound on the constraint functions. The approach in \cite{singh2025stochastic} is only applicable to the deterministic version of \eqref{mainProb} and achieves suboptimal $\O{\frac{1}{\eps^2}}$ complexity, unless the dual iterates are assumed bounded. 

Primal algorithms inspired from Polyak's subgradient iteration have been widely applied to solve \eqref{mainProb}  in~\cite{nedic2019random, necoara2022stochastic, singh2024unified, singh2024mini, singh2024stochastic, singh2024stochastic_movingball}. \col{Among these, \cite{necoara2022stochastic} is the closest to the present work, since it also considers stochastic composite optimization with functional constraints, uses only one sampled constraint per 
iteration, and attains $\O{\frac{1}{\eps}}$ SFO complexity in the  strongly convex setting.} Broadly, these approaches have two key steps: (a) a standard (stochastic) subgradient step and (b) a projection of  this step onto a halfspace approximation of a random constraint  function. \col{The proposed HPS framework is related in spirit, since it also uses only one sampled constraint per iteration, but applies a proximal step with respect to a linearized hinge penalty together with the nonsmooth term $h$. Thus, while \cite{necoara2022stochastic} performs a Polyak-type feasibility correction for the sampled constraint, HPS performs a stochastic proximal update for an exact hinge-penalty reformulation. Interestingly, this seemingly small modification allows us to use smooth non-Lipschitz constraint functions while also yielding improved rates in the finite-sum setting.}

\col{The assumption regimes in the two works are also different. The analysis in \cite{necoara2022stochastic} requires the functional constraints to be globally Lipschitz continuous, or equivalently to have uniformly bounded subgradients, which is natural in nonsmooth settings but excludes smooth constraints whose gradients grow with $\norm{\x}$. In contrast, HPS, VR-HPS, and N-HPS are designed for smooth constraint functions, such as squared-error and quadratic constraints, which are generally not globally Lipschitz continuous. Also, the $\O{\frac{1}{\eps}}$ rate in \cite{necoara2022stochastic}  relies on the linear-regularity condition. The baseline HPS achieves the same dependence on $\epsilon$ under Slater's condition alone, although with a worse dependence on $m$. The proposed N-HPS also achieves almost the same rate for smooth non-Lipschitz constraint functions under linear regularity, thereby matching the SFO complexity of \cite{necoara2022stochastic} while relaxing the global Lipschitz constraint assumption. For the finite-sum case, the unified variance reduction framework in \cite{singh2024unified} achieves only suboptimal $\O{\frac{1}{\eps}}$ rates. In comparison, the proposed VR-HPS achieves an SFO complexity of $\O{\frac{1}{\sqrt{\eps}}}$ for smooth objective and constraints.}


Exact penalty methods and specifically the hinge penalty function has been extensively studied \cite[Sec. 4.3.1]{bertsekas}, \cite{han1979exact}. Various approaches exist for solving the resulting problems, including the model-based method \cite{davis2019stochastic}, prox-linear method \cite{zhang2021stochastic}, smooth penalty method \cite{xiao2019penalized,thomdapu2019optimal}, and stochastic sequential quadratic programming framework~\cite{sanyal2025stochastic}. The present work is also related to these, as the HPS algorithm begins with reformulating \eqref{mainProb} using the exact penalty method and subsequently utilizes a stochastic linearized proximal method to solve it, while making only one SFO call per-iteration. The subsequent algorithms VR-HPS and N-HPS deviate further from this template, leading to novel algorithms and different rates. 

\begin{table}[t]
\centering
\footnotesize
\caption{\small Algorithms for \eqref{mainProb} with stochastic
objective; big-$\tilde{\mathcal{O}}$ notation hides logarithmic terms. The dependence on $m^2$ for \cite{xu2020primal} is not explicitly stated but is derived for \eqref{mainProb} in Appendix-\ref{appendix: SFO}.}
\label{tab:sc-constraints1}
\setlength{\tabcolsep}{3pt}
\renewcommand{\arraystretch}{1.5}
\begin{tabularx}{\columnwidth}{@{}l p{1.9cm} p{1.9cm}  p{1.9cm} p{2.5cm}@{}}
\toprule
Ref. & Regularity Assumption & Assumption \newline on $g_j$ &  Memory \newline Complexity & SFO \newline Complexity \\
\midrule
\cite{xu2020primal} & \hspace{4mm}$\times$ & Lipschitz & $\O{m}$ & \col{$\Ot{m^2\epsilon^{-1}}$} \\
\cite{necoara2022stochastic} & \hspace{4mm}$\checkmark$ & Lipschitz & $\O{1}$ & $\O{\epsilon^{-1}}$ \\
HPS  & \hspace{4mm}$\times$ & Smooth & $\O{1}$ & \col{$\O{m^2\epsilon^{-1}}$} \\
N-HPS  & \hspace{4mm}$\checkmark$ & Smooth & $\O{1}$ & $\Ot{\eps^{-1}}$ \\
\bottomrule
\end{tabularx}
\end{table}

\begin{table}[t]
\centering
\footnotesize
\caption{\small Algorithms for \eqref{mainProb} with deterministic or finite-sum objective; both smoothness and Lipschitz
continuity are required in \cite{singh2025stochastic,
singh2024stochastic_movingball}; better rates are obtained in \cite{singh2025stochastic} if the dual iterates are assumed bounded. The dependence on $n$ and $m$ for some algorithms is not explicitly stated but is derived for \eqref{mainProb} in in Appendix-\ref{appendix: SFO}.}
\label{tab:sc-constraints2}
\setlength{\tabcolsep}{3pt}
\renewcommand{\arraystretch}{1.5}
\begin{tabularx}{\columnwidth}{@{}l p{1.2cm} p{1cm} p{1.3cm} p{1.2cm} p{2.5cm}@{}}
\toprule
Ref. & Regularity Assump. & Finite-Sum & Assump. \newline on $g_j$ & Memory \newline Complex. & SFO \newline Complex.\\
\midrule
\cite{singh2025stochastic} & \hspace{4mm}$\times$ & \hspace{2mm}$\times$ & Lipschitz & $\mathcal{O}(m)$ & \col{$\O{nm^2\epsilon^{-2}}$} \\
\cite{nedic2019random} & \hspace{4mm}$\checkmark$ & \hspace{2mm}$\times$ & Lipschitz & $\O{1}$ & \col{$\O{n\epsilon^{-1}}$} \\
\cite{singh2024stochastic} & \hspace{4mm}$\checkmark$ & \hspace{2mm}$\times$ & Lipschitz & $\O{1}$ & \col{$\O{n\epsilon^{-1}}$} \\
\cite{singh2024stochastic_movingball} & \hspace{4mm}$\checkmark$ &  \hspace{2mm}$\times$ & Lipschitz & $\O{1}$ & \col{$\O{n\epsilon^{-1}}$} \\
\cite{singh2024mini} & \hspace{4mm}$\checkmark$ & \hspace{2mm}$\checkmark$ & Lipschitz & $\O{1}$ & \col{$\O{(n+m)\epsilon^{-1}}$} \\
\cite{singh2024unified} & \hspace{4mm}$\checkmark$ & \hspace{2mm}$\checkmark$ & Lipschitz & $\O{1}$ & \col{$\O{n+\epsilon^{-1}}$} \\
VR-HPS & \hspace{4mm}$\times$ & \hspace{2mm}$\checkmark$ & Smooth & $\mathcal{O}(m)$ & \col{$\O{n+\frac{m+\sqrt{n}}{\sqrt{\epsilon}}}$} \\
\bottomrule
\end{tabularx}
\end{table}

\subsection{Notations and Organization}
Regular (non-bold) letters denote scalars, while bold-faced letters represent column vectors. For a scalar value $v$, we define  $[v]_+:= \max\{0,v\}$, The Euclidean norm of a vector $\x$ is denoted as $\norm{\x}$. The proximal operator with respect to a function $h$ is defined as
\begin{align}
    \prox{h}{\z} := \arg\min_{\u} \frac{1}{2}\norm{\z-\u}^2 +  h(\u).
\end{align}
We use the notation $\nabla h(x)$ to denote an arbitrarily chosen subgradient in $\partial h(x)$ whenever $h$ is non-smooth.

The remainder of this paper is organized as follows: Sec. \ref{sec:proposed_algos} details the proposed HPS algorithm, Sec. \ref{sec:vrhps} discusses VR-HPS algorithm, and Sec. \ref{sec:n-hps} introduces N-HPS algorithm. The robust regression problem is introduced and all three algorithms are evaluated against other benchmarks in Sec. \ref{sec-numerical}. Finally, Sec. \ref{sec-conclusion} concludes the paper. 

\section{Background and Assumptions} 
In this section, we introduce the exact penalty framework that will serve as the main vehicle for designing our algorithms. We also state the various assumptions required for the proposed algorithms. 

The exact penalty method for constrained optimization has been widely studied and incorporates the constraints directly into the objective via a nonsmooth penalty function \cite[Sec. 4.3.1]{bertsekas}. Using this approach, we reformulate \eqref{mainProb} as
\begin{align}
	\min_{\x\in \Rn^d} ~ F(\x) := f(\x) + h(\x) + \frac{\gamma}{m} \sum_{j = 1}^{m}[g_j(\x)]_+
	\tag{$\mathcal{P}_1$}\label{penalty1}
\end{align}
where $\gamma > 0$ is the penalty parameter. Observe that instead of enforcing the constraints via projections onto the feasible region, the penalty reformulation attempts to correct the constraint violations through hinge penalties applied to individual constraints. Remarkably, and unlike various other smooth penalty functions, the exact penalty ensures that any solution to \eqref{penalty1} is also a solution to \eqref{mainProb}, provided that $\gamma$ is sufficiently large \cite[Prop. 4.3.2]{bertsekas}. We remark that other types of penalty functions, such as $\max_k [g_k(\x)]_+$, have also been used in the literature \cite{bertsekas1975necessary, sanyal2025stochastic}. Here, we utilize a separate penalty for each constraint in \eqref{penalty1}, resulting in a finite-sum structure that will be the key to developing stochastic algorithms capable of handling a large number of constraints. 

For brevity, we collect the non-smooth terms into 
\begin{align}
    \phi(\x) &\coloneqq \tfrac{1}{m}\sum_{j=1}^m \phi_j(\x) & \phi_j(\x) := h(\x) + \gamma[g_j(\x)]_+,
\end{align}
for $1\leq j \leq m$ so that $F(\x)=f(\x) + \phi(\x)$. The somewhat redundant appearance of $h$ in each $\phi_j$ is intentional, as it allows us to work with a collection of per-constraint non-smooth components $\{\phi_j\}_{j=1}^m$ and apply stochastic algorithms that access only one $\phi_j$ per-iteration.

The feasible region for the $j$-th constraint is denoted by $\cX_j := \{\x \mid g_j(\x) \leq 0\}$ and the entire feasible region is $\cX = \cap_{j=1}^m \cX_j$. A point $\x_\star$ is optimal if it satisfies the first order optimality condition
\begin{align}
      \ip{\nabla f(\x_\star) + \nabla h(\x_\star)}{\tilde{\x}-\x_\star} &\geq 0 & \forall ~\tilde{\x}\in \cX  \label{optgrad0}
\end{align}
for some $\nabla h(\x_\star)\in \partial h(\x_\star)$. Since we consider $f$ to be strongly convex, the optimal $\x_\star$ is guaranteed to be unique. For stochastic optimization algorithms, a random point $\x$ is near-optimal if its mean squared distance from $\x_\star$ is small, i.e., $\EE\norm{\x - \x_\star}^2 \leq \epsilon$. 

We assume access to a stochastic first order (SFO) oracle that returns $\{\nabla f_{i_t}(\x), \nabla g_{j_t}(\x),  g_{j_t}(\x)\}$ for a given $\x$ and random indices $i_t$ and $j_t$. The SFO complexity of various algorithms will be measured by the total number of oracle calls required to obtain a near-optimal point. 

\subsection{Assumptions} \label{sec:assumption}
We first discuss assumptions that are common to all the algorithms. We begin with the standard assumption regarding the smoothness and convexity of the various functions. 
\begin{assumption}\label{a1}
	The functions $\{f_i\}_{i=1}^n$ and $\{g_j\}_{j=1}^m$ are proper, closed, convex, and $L_f$- and $L_g$-smooth, respectively. The regularizer $h$ is \col{non-negative,} proper, closed, and convex, but may be non-smooth.
\end{assumption}
Assumption \ref{a1} delineates the structural properties of the problem class under consideration. Note that we do not require $h$ to be Lipschitz and hence allow $h$ to contain indicator functions of convex sets. Next, we require the Slater condition which is one of the simplest and most widely used constraint qualifications (CQ) for convex optimization problems. 
\begin{assumption}\label{slater}
	Slater condition holds for \eqref{mainProb}, i.e., there exists a feasible $\tx$ such that 
	\begin{align}
		g_j(\tx) &\leq -\nu < 0& 1\leq j \leq m	
	\end{align}
	Additionally, we assume that the optimality gap at the Slater point is bounded, i.e., $f(\tx) + h(\tx) - f(\x_\star) - h(\x_\star) \leq \bt$.
\end{assumption}
Assumption \ref{slater} implies the existence of a primal-dual optimum pair $(\x_\star,\lam_\star)$ with $\norm{\lam_\star}_1 \leq \frac{\bt}{\nu}$; see \cite{sanyal2025stochastic}. Further, as shown in Appendix \ref{penaltyequivalent}, choosing $\gamma > m\frac{\bt}{\nu}$ ensures that a KKT point of \eqref{mainProb} is also a KKT point of \eqref{penalty1} and hence $\x_\star = \arg\min_{\x} F(\x)$. Under Assumption \ref{slater}, the optimality condition in \eqref{optgrad0} becomes $0 \in \partial F(\x_\star)$, or equivalently, 
\begin{align}
	\nabla f(\x_\star) + \nabla \phi(\x_\star) = 0. \label{optgrad2}
\end{align}
for some $\nabla \phi(\x_\star) \in \partial \phi(\x_\star)$. Next, we have the strong convexity assumption. 
\begin{assumption}\label{sc}
	The objective function $f(\x)$ is $\mu$-strongly convex.
\end{assumption}
In this context, we define the objective condition number $\kappa_f := \frac{L_f}{\mu}$ and the constraint condition number $\kappa_g := \frac{L_g}{\mu}$. Note that for problems satisfying Assumption \ref{sc}, the optimal $\x_\star$ is unique. Combining Assumption \ref{sc} with the bound in Assumption \ref{slater} yields
\begin{align}
    \frac{\mu}{2}\norm{\tx - \x_\star}^2 \leq f(\tx) + h(\tx) - f(\x_\star) - h(\x_\star) \leq \bt. \label{optdist}
\end{align}

For vanilla SGD variants, we also need the following bounded variance assumption. 
\begin{assumption}\label{sig}
	The stochastic gradients $\{\nabla f_i(\x_\star)\}_{i=1}^n$ have bounded variance at $\x_\star$, i.e., 
	\begin{align}
		\En{\norm{\nabla f_{i_t}(\x_\star) - \nabla f(\x_\star)}^2}\leq \sigma_\star^2	\label{sigeq}
	\end{align}
	 where the expectation is with respect to the random index $i_t$.  
\end{assumption}
The variance-reduced SGD proposed in Sec. \ref{sec:vrhps} will however not require Assumption \ref{sig}. For some of the algorithms, we need the gradients of the constraint functions to be bounded at $\x_\star$. 
\begin{assumption}\label{gstarbound}
	The constraint function gradients are mean-square bounded at $\x_\star$, i.e., $\En{\norm{\nabla g_{j_t}(\x_\star)}^2}\leq G_\star^2$. Likewise, the function gradients are bounded at $\x_\star$, i.e., $\norm{\nabla f(\x_\star)}^2 + \norm{\nabla h(\x_\star)}^2 \leq B_{\star}$ for some $\nabla h(\x_\star) \in \partial h(\x_\star)$. 
\end{assumption} 

Finally, Sec. \ref{sec:n-hps} will require a stronger linear regularity assumption, also considered in \cite{patrascu2021stochastic}. 
Let $\pi_{\cX}(\x)$ and $\pi_{\cX_{j_t}}(\x)$ be the projections of $\x$ onto the sets $\cX:=\{\u | g_j(\u)\leq 0, 1\leq j \leq m\}$ and $\cX_{j_t}:=\{\u | g_{j_t}(\u)\leq 0\}$, respectively. 
\begin{assumption}\label{regular}
	The constraints satisfy linear regularity with parameter $r$, i.e., $r\norm{\x - \pi_{\cX}(\x)}^2 \leq \Ej{\norm{\x - \pi_{\cX_{j_t}}(\x)}^2}$, where $\Ej{\cdot}$ denotes the expectation with respect to the random index $j_t$. 
\end{assumption}
Assumption \ref{regular} states that the expected distance to $\cX$ is controlled by the expectation of distances to the individual sets $\cX_{j_t}$, thereby excluding pathological cases involving nearly parallel or redundant constraints. \col{While Slater's CQ implies a bounded-set version of 
linear regularity through classical error-bound results \cite{bauschke1999strong}, it yields constants that depend on the bounded region. In particular, for any bounded set $\cS$, one obtains a constant $r_{\cS} > 0$ such that $r_{\cS}\norm{\x - \pi_{\cX}(\x)}^2 \leq \Ej{\norm{\x - \pi_{\cX_{j_t}}(\x)}^2}$ for all $\x \in \cS$. In contrast, the global linear regularity assumed in  Assumption~\ref{regular} and in \cite{wang2016stochastic2} requires the inequality to hold for the entire, possibly unbounded, decision region with a uniform constant, and is therefore stronger. In particular, such a global constant cannot be bounded using only the Slater parameters $\nu$ and $\bt$. A related residual-based global regularity condition $\norm{\x - \pi_{\cX}(\x)}^2 \leq c\Ej{([g_{j_t}(\x)]_+)^2}$ is also assumed in \cite{nedic2019random, necoara2022stochastic, singh2024unified, singh2024mini, singh2024stochastic, singh2024stochastic_movingball} and it also does not follow from Slater's CQ alone.}


We conclude this section by mentioning that the various preliminary inequalities that follow from Assumptions \ref{a1}-\ref{regular} are listed in Appendix \ref{basic} and will be repeatedly used.

\section{Hinge-Proximal SGD} \label{sec:proposed_algos}
In this section, we introduce the hinge-proximal framework for solving \eqref{penalty1}, which exploits the finite-sum structure of $\phi$ to access only one constraint per-iteration, while still operating in the primal domain. Our starting point is the stochastic splitting proximal gradient (SSPG) algorithm of \cite{Patrascu2021}, which also handles stochastic composite problems such as \eqref{penalty1}. At iteration $t$, the SSPG algorithm for solving \eqref{penalty1} would involve (a) performing an SGD update to calculate $\z_t = \x_t - \eta_t\nabla f_{i_t}(\x_t)$ for a random index $i_t$; and then (b) calculating $\x_{t+1} = \prox{\eta_t \phi_{j_t}}{\z_t}$. 

For SSPG to be efficient, $\phi_{j_t}$ must be proximally tractable, which is often not the case in many instances of \eqref{mainProb}. The proposed HPS algorithm retains the same two-stage structure, but replaces the difficult proximal step with a simpler proximal with respect to a partially linearized penalty $\tphi_j(\u,\x_t) \coloneqq h(\u) + \gamma [g_j(\x_t) + \ip{\nabla g_j(\x_t)}{\u-\x_t}]_+$, resulting in the update
\begin{align} 
	\x_{t + 1} = \prox{\eta_t\tphi_{j_t}}{\z_t} &= \argmin_{\u} \tfrac{1}{2\eta_t}\norm{\z_t - \u}^2 \label{xupdate-sgd}\\
    &\hspace{-1cm}+ h(\u) + \gamma [g_{j_t}(\x_t) + \ip{\nabla g_{j_t}(\x_t)}{\u-\x_t}]_+. \nonumber
\end{align}
Since \eqref{xupdate-sgd} involves the use of proximal step with respect to the sum of $h$ and a hinge function, we refer to it as hinge-proximal SGD. Observe that $\tphi_{j_t}$ retains $h$ as is, but penalizes violations of the sampled constraint $g_{j_t}$ as measured by its linear proxy. Unlike SSPG, the HPS updates do not require $g_{j_t}$ to be proximally tractable and only use $\nabla g_{j_t}(\x_t)$. We show in the Appendix \ref{tractable} that, whenever $h$ is proximally tractable, $\prox{\eta_t \tilde \phi_{j_t}}{\cdot}$ can be evaluated to accuracy $\epsilon$ using $O(\log(1/\epsilon))$ calls to the proximal operator of $h$. In the special case when $h = 0$, the update in \eqref{xupdate-sgd} reduces to projection onto a half-space, which can be written in closed-form. The full algorithm is summarized in Algorithm \ref{algo:HPS}. 

\begin{algorithm}
	\caption{Hinge-Proximal SGD}
	\label{algo:HPS}
	\begin{algorithmic}[1]
		\State \textbf{Parameters:} Step-size schedule $\eta_t$, penalty factor $\gamma$
		\State \textbf{Input:} initial point $\x_1$
		\For {$t = 1, \ldots, T-1$ }
		\State Evaluate $\{\nabla f_{i_t}(\x_t), \nabla g_{j_t}(\x_t)\}$
        \State Set $\z_t \gets \x_t - \eta_t\nabla f_{i_t}(\x_t)$
		\State Update $\x_{t + 1} \gets \prox{\eta_t\tphi_{j_t}}{\z_t}$ \label{mainAlgostep-impr}
		\EndFor
		\State \textbf{Output:} $\x_T$
	\end{algorithmic}
\end{algorithm}

\col{While Algorithm~\ref{algo:HPS} is motivated by SSPG~\cite{Patrascu2021},
the substitution of $\prox{\eta_t\phi_{j_t}}{\cdot}$ by $\prox{\eta_t\tphi_{j_t}}{\cdot}$
is not a routine modification. The standard three-point property for
proximal operators, which underlies the SSPG analysis, does not carry
over to the linearized penalty $\tphi_{j_t}$. Lemma~\ref{lem0}
establishes the hinge-proximal analogue that makes the convergence
analysis tractable, and is the key technical departure from prior
proximal splitting approaches~\cite{Patrascu2021}.}

\col{HPS can also be interpreted as a stochastic version of the prox-linear method \cite{davis2019stochastic}, applied to the exact-penalty objective $F(\x)=f(\x)+h(\x)+\frac{\gamma}{m}\sum_{j=1}^m[g_j(\x)]_+$. Indeed, linearizing the sampled smooth objective $f_{i_t}$ and the sampled hinge term $[g_{j_t}]_+$ at $\x_t$, while keeping $h$ proximal, yields  the update in Algorithm~\ref{algo:HPS}. The key difference from a deterministic prox-linear step on $F$, as in \cite{davis2019stochastic}, is that HPS accesses only one constraint term $g_{j_t}$ per iteration, rather than linearizing all $m$ hinge terms.}

To analyze Algorithm~\ref{algo:HPS}, we first state a key per-iteration inequality for a generic hinge-proximal update. We deliberately present it in a more general form, in terms of an auxiliary point $\w_t$, so that it can be reused in the analysis of other algorithms in later sections. The inequality \eqref{eqlem0}  can be viewed as a hinge-proximal analogue of the classical three-point property for proximal mappings, and will serve as a common building block for the basic HPS recursion and its variance-reduced version. 
\begin{lemma}\label{lem0}
	Under Assumption \ref{a1}, the update $\w_t = \prox{\eta_t\tphi_{j_t}}{\z_t}$ satisfies
	\begin{align}\label{eqlem0}
		&\norm{\w_t - \x_\star}^2 \leq \norm{\x_t - \x_\star}^2 - (1 - \eta_t \gamma L_g)\norm{\w_t - \x_t}^2  \\ 
		& \quad  + 2\ip{\x_t-\z_t}{\x_\star - \w_t} + 2\eta_t (\phi_{j_t}(\x_\star) - \phi_{j_t}(\w_t)). \nonumber
	\end{align}
\end{lemma}
\begin{IEEEproof}
	Since the objective in \eqref{xupdate-sgd} is $\tfrac{1}{\eta_t}$-strongly convex, we have 
	\begin{align} \label{three-point-ineq1}
		& \tfrac{1}{2\eta_t}\norm{\w_t - \z_t}^2 + \tphi_{j_t}(\w_t,\x_t) + \tfrac{1}{2\eta_t}\norm{\w_t - \x_\star}^2\nonumber\\
		& \leq \tfrac{1}{2\eta_t}\norm{\x_\star - \z_t}^2 + \tphi_{j_t}(\x_\star,\x_t).
	\end{align}
	Multiplying by $2\eta_t$, inserting $\x_t$ into the norm terms, and canceling the common terms from both sides, we obtain
	\begin{align}
		&\norm{\w_t - \x_t}^2 + 2\eta_t\tphi_{j_t}(\w_t,\x_t) + \norm{\w_t - \x_\star}^2\nonumber\\
		& \leq \norm{\x_\star - \x_t}^2 + 2\ip{\x_t-\z_t}{\x_\star-\w_t} + 2\eta_t\tphi_{j_t}(\x_\star,\x_t) \nonumber\\
		&\leq \norm{\x_\star - \x_t}^2 + 2\ip{\x_t-\z_t}{\x_\star-\w_t} + 2\eta_t\phi_{j_t}(\x_\star)\label{three-point-ineq2}
	\end{align}
	where we have used the convexity of $g_j$ which implies that $\tg_j(\x_\star,\x_t) \leq g_j(\x_\star)$. Finally, using the $L_g$-smoothness of $g_{j_t}$ (Assumption~\ref{a1}) and the standard quadratic upper bound, we obtain $[g_{j_t}(\w_t)]_+ \leq [\tg_{j_t}(\w_t,\x_t)]_++\frac{L_g}{2}\norm{\w_t - \x_t}^2$. Substituting this into \eqref{three-point-ineq2} and rearranging, we obtain the required bound. 
\end{IEEEproof}

We expect this inequality to be useful beyond the present setting, for example in analyzing other stochastic algorithms that alternate between smooth updates and hinge-type corrections. Having established the preliminary lemma, we are now ready to establish the SFO complexity of Algorithm \ref{algo:HPS}.

\begin{theorem}\label{hps-thm} 
	Under Assumptions \ref{a1},\ref{slater},\ref{sc},\ref{sig}, and \ref{gstarbound}, for \col{$\eta_t =  \frac{2}{\mu t+2\tilde{L}}$ where $\tilde{L} = 2\max\{\gamma L_g, 4L_f\}$} and $\gamma = \frac{2m\tilde{B}}{\nu}$, and for $\norm{\x_1-\x_\star}^2 \leq B_x$, the output of Algorithm \ref{algo:HPS} satisfies: 
	\begin{align}
		\E{\norm{\x_T - \x_\star}^2} &= \O{\frac{G_\star^2 m^2 + \sigma_\star^2}{\mu^2 T}+\frac{\kappa_g^2 m^2+\kappa_f^2}{T^2}}
	\end{align}
	resulting in the SFO complexity of $\O{\frac{G_\star^2m^2 +\sigma_\star^2}{\mu^2\epsilon} + \frac{m\kappa_g + \kappa_f}{\sqrt{\epsilon}}}$. 
\end{theorem}
\begin{IEEEproof}
    The proof involves plugging in the HPS update into Lemma~\ref{lem0}, bounding the resulting terms using Assumptions~\ref{sig}–\ref{gstarbound}, and then using the resulting recursion to obtain the final rate. 

	Writing \eqref{eqlem0} for $\z_t = \x_t  - \eta_t \nabla f_{i_t}(\x_t)$ and $\w_t = \x_{t+1}$, we obtain
	\begin{align}\label{hps-proof}
		&\norm{\x_{t+1} - \x_\star}^2 \leq \norm{\x_t - \x_\star}^2 - (1 - \eta_t \gamma L_g)\norm{\x_{t+1} - \x_t}^2  \nonumber\\ 
		&  + 2\eta_t \ip{\nabla f_{i_t}(\x_t)}{\x_\star - \x_t} + 2\eta_t (\phi_{j_t}(\x_\star) - \phi_{j_t}(\x_{t+1})) \nonumber\\
		&  + 2\eta_t \ip{\nabla f_{i_t}(\x_t)}{\x_t - \x_{t+1}}.
	\end{align}
	Firstly, from the convexity of $\phi_{j_t}$ and the fact that $j_t$ is independent of $\x_t$, we have that 
	\begin{align}\label{hps-proof2}
		&\E{\phi_{j_t}(\x_\star) - \phi_{j_t}(\x_{t+1})} \leq \E{\ip{\nabla \phi_{j_t}(\x_\star)}{\x_\star - \x_{t+1}}} \\
		&\eqtext{\eqref{optgrad2}} -\E{\ip{\nabla f(\x_\star)}{\x_\star - \x_t }} + \E{\ip{\nabla \phi_{j_t}(\x_\star)}{\x_t - \x_{t+1}}}\nonumber
	\end{align}
	Taking expectation in \eqref{hps-proof} and substituting \eqref{hps-proof2}, we obtain 
	\begin{align}
		&\EE\norm{\x_{t+1} - \x_\star}^2 \leq \En{\norm{\x_t - \x_\star}^2} \nonumber\\
		&-(1 - \eta_t \gamma L_g)\En{\norm{\x_{t+1} - \x_t}^2} \nonumber\\
		&-2\eta_t \E{\ip{\nabla f_{i_t}(\x_\star) - \nabla f_{i_t}(\x_t)}{\x_\star - \x_t}} \nonumber\\
		&+ 2\eta_t \E{\ip{\nabla f_{i_t}(\x_t) + \nabla \phi_{j_t}(\x_\star)}{\x_t - \x_{t+1}}}.\label{hps-proof3}
	\end{align}
	Let us now bound the last two terms separately. \col{The first term can be bounded using Assumption \ref{sc}. Since $i_t$ is independent of $\x_t$, we have the averaged relation
    \begin{align}
        &\E{\ip{\nabla f_{i_t}(\x_\star) - \nabla f_{i_t}(\x_t)}{\x_\star - \x_t}} \nonumber\\
        &= \E{\ip{\nabla f(\x_\star) - \nabla f(\x_t)}{\x_\star - \x_t}} \\
        &= \EE D_f(\x_t,\x_\star) + \EE D_f(\x_\star,\x_t) \\
        &\geqtext{\eqref{bound_sc}}  \EE D_f(\x_t,\x_\star) + \tfrac{\mu}{2} \EE\norm{\x_t -\x_\star}^2\label{hps-t1}
    \end{align}
    where $D_f(\x,\y) := f(\x) - f(\y) - \ip{\nabla f(\y)}{\x-\y}$ is the Bregman divergence between $\x$ and $\y$ with respect to $f$. The last term in \eqref{hps-proof3} can be bounded using \eqref{young} as
    \begin{align}
        & \EE\ip{\nabla f_{i_t}(\x_t) + \nabla \phi_{j_t}(\x_\star)}{\x_t - \x_{t+1}} \leqtext{\eqref{young}} \tfrac{1}{4\eta_t}\EE\norm{\x_t - \x_{t+1}}^2 \nonumber\\
		& + 2\eta_t\EE\big[\norm{\nabla f_{i_t}(\x_t)-\nabla f(\x_\star)}^2 + \norm{\nabla \phi_{j_t}(\x_\star) - \nabla \phi(\x_\star)}^2\big] \nonumber\\
        &\leqtext{\eqref{fvar},\eqref{phistarbound}} \tfrac{1}{4\eta_t}\EE\norm{\x_t - \x_{t+1}}^2 \nonumber\\
        & \hspace{1cm}+ 2\eta_t(4 L_f\EE D_f(\x_t,\x_\star) + 2\sigma_\star^2 + \gamma^2G_\star^2).\label{hps-t2}
    \end{align}
    Substituting \eqref{hps-t1}-\eqref{hps-t2} into \eqref{hps-proof3}, we obtain
    \begin{align}
		&\E{\norm{\x_{t+1} - \x_\star}^2} \leq \left(1-\mu \eta_t \right)\En{\norm{\x_t - \x_\star}^2} \nonumber\\
        &+ 4\eta_t^2(\gamma^2G_\star^2 + 2\sigma_\star^2) -(\tfrac{1}{2} - \eta_t \gamma L_g)\En{\norm{\x_{t+1} - \x_t}^2} \nonumber\\
		&- 2\eta_t\left(1 - 8\eta_t L_f\right)\EE D_f(\x_t,\x_\star)\label{hps-proof4}
	\end{align}
    For $\eta_t \leq \frac{1}{\tilde{L}}$ where $\tilde{L} = 2\max\{\gamma L_g, 4L_f\}$, both coefficients 
$(\tfrac{1}{2}-\eta_t\gamma L_g)$ and $(1-8\eta_tL_f)$ are nonnegative, so the last two terms in \eqref{hps-proof4} can be dropped. Hence, using \eqref{rect} (for $\sY_t = 0$), we obtain
	\begin{align}
		\EE\norm{\x_T - \x_\star}^2 \leq  16\tfrac{(\gamma^2G_\star^2 + 2\sigma_\star^2)}{\mu^2T} + 16\tfrac{\max\{\gamma^2L_g^2, 16L_f^2\} B_x}{\mu^2T^2} \label{hps-bxbound}
	\end{align}
    for initial $\x_1$ such that $\norm{\x_1 - \x_\star}^2\leq B_x$. } Again note that for $\x_\star$ to be the minimizer of \eqref{mainProb}, we need $\gamma  > \frac{m\tilde{B}}{\nu}$, which implies that the squared distance from the optimum is bounded as $\EE\norm{\x_T - \x_\star}^2 \leq \O{\frac{G_\star^2 m^2 + \sigma_\star^2}{\mu^2 T}+\frac{\kappa_g^2 m^2+\kappa_f^2}{T^2}}$
	where $\O{\cdot}$ hides the $\tilde{B}$, $B_x$, and $\nu$ terms, as well as the universal constants. The resulting iteration complexity is given by $\O{\frac{G_\star^2m^2 +\sigma_\star^2}{\mu^2\epsilon} + \frac{m\kappa_g + \kappa_f}{\sqrt{\epsilon}}}$
\end{IEEEproof}

\col{Note that the condition $\norm{\x_1-\x_\star}^2\le B_x$ in the statement of Theorem~\ref{hps-thm} is primarily used to quantify the effect of initialization on the iteration complexity. In particular, from \eqref{hps-bxbound}, $B_x$ appears only in the lower-order term and not in the leading $\O{1/T}$ term. Thus, the leading iteration complexity is insensitive to the choice of initialization. In practice, $\x_\star$ is always unknown, and for any chosen initialization, the bound holds for some finite problem-dependent constant $B_x$, which appears only in the convergence bounds and is never used by the algorithm.}

The result in Thm. \ref{hps-thm} also extends to nonsmooth $f$ under a Lipschitz subgradient assumption, as shown in in Appendix-\ref{appendix: HPS_nonsmooth_f}. Recall that the SFO complexity of SGD for solving unconstrained optimization problems is $\O{\tfrac{\sigma_\star^2}{\mu^2\epsilon}}$. Thus, when $m$ is treated as a constant, the SFO complexity of HPS matches that of SGD, up to constants, despite the presence of nonlinear functional constraints and using only one constraint per-iteration. In contrast to \cite{necoara2022stochastic}, which attains a similar rate under a linear regularity assumption, Thm. \ref{hps-thm} only requires Slater's condition. From Thm. \ref{hps-thm}, we can also see that $\E{\text{dist}(\x_T, \cX)^2} \leq \En{\norm{\x_T - \x_\star}^2} = \O{1/T}$, in line with the results obtained in \cite{necoara2022stochastic}. However, establishing an $\O{1/T}$ bound on the maximum constraint violation $\max_j [g_j(\x_T)]_+$ seems difficult in the current context and likely requires additional assumptions. \col{We also remark that the $\O{T^{-1/2}}$ result for the convex case does not follow from our current assumptions. While these assumptions can ensure that $\EE D_f(\hat{\x}_T,\x_\star)=\O{T^{-1/2}}$ for the averaged iterate $\hat{\x}_T$, they do not imply analogous bounds on the optimality gap or constraint violation. It appears that the convex case requires additional assumptions, such as bounded iterates, bounded penalty subgradients, or stochastic bounded gradients as in \cite{necoara2022stochastic}. }

The remaining drawback is the $\O{m^2}$ dependence inherited from the penalty parameter $\gamma$ and the variance terms involving $\sigma_\star$ and $G_\star$, which dominate the rate. Hence, Thm. \ref{hps-thm} highlights a clear target for improvement: in the finite-sum setting, applying variance reduction on both the smooth and nonsmooth components should remove the $m^2$ factor and improve over the suboptimal $\O{1/\eps}$ rate. We develop such a variance-reduced method next. 



\section{Variance Reduced Hinge-Proximal SGD}\label{sec:vrhps}
In this section, we construct a variance-reduced version of HPS for the finite-sum setting, i.e., for $f(\x) = \frac{1}{n}\sum_{i=1}^n f_{i}(\x)$. The proposed VR-HPS algorithm achieves an improved \col{$\O{n+\frac{m+\sqrt{n}}{\sqrt{\epsilon}}}$} complexity, which to the best of our knowledge, is the first such rate for solving \eqref{mainProb} while accessing only one constraint gradient per-iteration. Our design combines SVRG-style variance reduction for the smooth component $f$ with a decoupled variance reduction scheme for the hinge penalty $\phi$, adapted from \cite{mishchenko2019stochastic} but modified to work with linearized constraints and without proximal access to each $g_j$.

For the smooth term, we use a SVRG-type variance reduced estimator $\v_t \coloneqq \nabla f_{i_t}(\x_t) - \nabla f_{i_t}(\bx_t) + \nabla f(\bx_t)$, where $\bx_t$ is a checkpoint at which the full gradient $\nabla f(\bx_t)$ is periodically computed and kept fixed between updates. \col{For the non-smooth term, if the full penalty subgradient $\nabla \phi(\x_t)$ were available, the update would use the direction $\v_t+\nabla \phi(\x_t)$, leading to an intermediate point $\z_t=\x_t-\eta_t(\v_t+\nabla\phi(\x_t))$. However, computing $\nabla \phi(\x_t)=\frac{1}{m}\sum_{j=1}^m\nabla\phi_j(\x_t)$ requires accessing all the constraints at every iteration. Therefore, following the decoupling idea of \cite{mishchenko2019stochastic}, we introduce auxiliary variables $\y_{j,t}$ that track the individual penalty components and maintain their average $\y_t \coloneqq \tfrac{1}{m}\sum_{j=1}^{m}\y_{j,t}$ as an approximation of $\nabla\phi(\x_t)$.} At the $t$-th iteration, we only update $\y_{j_t,t+1}$ to replace $\y_{j_t,t}$ for a random $j_t$, so that $\y_t$ can be maintained as 
\begin{align}
    \y_{t + 1} &= \y_{t} + \frac{1}{m}(\y_{j_t, t + 1} - \y_{j_t, t}). \label{yupdate}
\end{align}
If $\y_t$ were the exact subgradient of $\phi(\x_t)$, the classical stochastic (sub-)gradient update would have involved updating $\x_t$ with $\x_t -\eta_t \v_t -  \eta_t \y_t$. Since $\y_t$ is only an approximate subgradient, additional correction terms need to be included. The update in \cite[Sec. 4]{mishchenko2019stochastic} is motivated by the connection between the subgradient and the proximal operation, and requires each regularizer to be proximally tractable. Here, we use a similar connection, but introduce a key innovation: we replace $g_{j_t}$ with its linearized version $\tg_{j_t}(\u,\x_t) := g_{j_t}(\x_t) + \ip{\nabla g_{j_t}(\x_t)}{\u-\x_t}$ within the update. Hence replacing $\phi_{j_t}$ by $\tphi_{j_t}(\u,\x_t)=h(\u) + \gamma [\tg_{j_t}(\u,\x_t)]_+$, we propose the update 
\begin{align}
    \x_{t + 1} &= \prox{\eta_t\tphi_{j_t}}{\x_t -\eta_t \v_t- \eta_t\y_t +\eta_t\y_{j_t,t}}\label{eq-vrhps-up1} \\
    &= \argmin_{\u\in\Rn^d} h(\u) + \gamma [g_{j_t}(\x_t) + \ip{\nabla g_{j_t}(\x_t)}{\u-\x_t}]_+ \nonumber\\
    & \hspace{1cm}+ \tfrac{1}{2\eta_t}\norm{\x_t -\eta_t \v_t - \eta_t\y_t +\eta_t\y_{j_t,t} - \u}^2.  \label{eq-vrhps-up2}
\end{align}
As explained earlier, such a hinge-proximal update can be easily carried out if $h$ is proximally tractable. Specifically, we can use a one-dimensional search method to solve \eqref{eq-vrhps-up1} to accuracy $\epsilon$ using at most $O(\log(1/\epsilon))$ proximal evaluations of $h$. Interestingly, the linearization step also makes the update for $\y_{j_t,t+1}$  different from that in \cite{mishchenko2019stochastic}, which is now written as
\begin{align}\label{yjtupdate}
    \y_{j_t, t + 1} = \y_{j_t, t} + \frac{1}{2\eta_t}(\x_t - \x_{t + 1}) - (\v_t + \y_t).
\end{align}
\col{This is the key point at which VR-HPS differs from the decoupling method of \cite{mishchenko2019stochastic}. That method assumes proximal access to each nonsmooth component, whereas here $\phi_j(\x)=h(\x)+\gamma[g_j(\x)]_+$ is not proximally tractable for a general smooth nonlinear $g_j$. We therefore update $\x_{t+1}$ using the linearized surrogate $\tphi_{j_t}(\cdot,\x_t)$, and $\y_{j_t,t+1}$ tracks the corresponding linearized proximal residual. This modification leads to the update in \eqref{yjtupdate}; in particular, the factor of $2$ in the denominator is absent in \cite[Alg.~1]{mishchenko2019stochastic} and is needed for the proof of Theorem \ref{vrhps-thm}. } Combining these ingredients, we obtain the SVRG-based VR-HPS method summarized in Algorithm~\ref{algo:VRHPS}. Observe that compared to Algorithm \ref{algo:HPS}, VR-HPS requires an additional $\O{m}$ storage to maintain $\{\y_{j,t}\}_{j=1}^m$. 	

\col{This also distinguishes VR-HPS from~\cite{singh2024unified}, where variance reduction is applied to the finite-sum objective but not to the functional-constraint terms. In contrast, VR-HPS maintains the auxiliary variables $\{\y_{j,t}\}_{j=1}^m$ to track and variance-reduce the hinge-penalty components, yielding the improved finite-sum complexity at the cost of additional $\O{md}$ storage.}

\begin{algorithm}
	\caption{VR-HPS}
	\label{algo:VRHPS}
	\begin{algorithmic}[1]
		\State \textbf{Parameters:} Step-size schedule $\eta_t$, penalty factor $\gamma$
		\State \textbf{Input:} initial point $\x_1=\bx_1$ 
		\State \textbf{Input:} $\y_{j, 1} = \y_1$ such that $\norm{\y_1-\y_{j,\star}}^2\leq B_y$ for all $j$
		\State Compute $\nabla f(\bx_1)$
		\For {$t = 1, \ldots, T-1$ }
		\State Sample $i_t$ and $j_t$ uniformly at random
		\State Update $\v_t = \nabla f_{i_t}(\x_t) - \nabla f_{i_t}(\bx_t) + \nabla f(\bx_t)$
		\State Update $\bx_{t+1} = \x_t$ with probability $\tfrac{1}{n}$, otherwise set $\bx_{t+1} = \bx_t$
		\State If $\bx_{t+1}$ has been updated, compute $\nabla f(\bx_{t+1})$
		\State Update $\x_{t + 1}$ using \eqref{eq-vrhps-up1} \label{mainAlgobstep}
		\State Update $\y_{j_t,t+1}$ as per \eqref{yjtupdate}
		\State Update $\y_t$ as per \eqref{yupdate}
		\EndFor
		\State \textbf{Output} $\x_T$
	\end{algorithmic}
\end{algorithm}

Matching the notation in Algorithm \ref{algo:VRHPS}, let us denote $\y_{j,\star} \in \partial \phi_j(\x_\star)$ and $\y_\star = \frac{1}{m}\sum_{j=1}^m \y_{j,\star}$, so that \eqref{optgrad2} can be written as $\nabla f(\x_\star) + \y_\star = 0$ for some choice of subgradients. The analysis of Algorithm \ref{algo:VRHPS} relies critically on the following lemma, which tracks how the subgradient estimates $\y_{j_t,t}$ approach the optimal subgradients $\y_{j_t,\star}$ with $t$. The proof of Lemma \ref{lem1} utilizes the update in \eqref{yjtupdate}, optimality condition \eqref{optgrad2}, and some algebraic manipulations, particularly the use of \eqref{young}, and is provided in Appendix \ref{prooflem1}.  
\begin{lemma}\label{lem1}
	Under Assumptions \ref{a1}-\ref{slater}, we have from \eqref{yjtupdate}:
	\begin{align}
		& 2\eta_t^2\sum_{j=1}^m\En{\norm{\y_{j, t + 1} - \y_{j, \star}}^2} - 2\eta_t^2\sum_{j=1}^m \En{\norm{\y_{j,t}-\y_{j,\star}}^2} \nonumber\\
		&\leq \tfrac{1}{2}\En{\norm{\x_t - \x_{t + 1}}^2}  + 2\eta_t^2\En{\norm{\v_t - \nabla f(\x_\star)}^2}  \nonumber\\
		&-  2\eta_t \E{\ip{\v_t + \y_t - \y_{j_t,t}  + \nabla \phi_{j_t}(\x_\star) }{\x_t-\x_{t+1}}}.\label{eqlem1}
	\end{align}
\end{lemma}
Having stated the preliminary lemma, we are now ready to prove the main result that establishes the $\O{1/T^2}$ bound on the squared distance to the optimum. 

\begin{theorem}\label{vrhps-thm} 
	Under Assumptions \ref{a1},\ref{slater}, and \ref{sc},  for \col{$\eta_t =  \frac{2}{\mu t+2\tilde{L}}$ where $\tilde{L} = 2\max\{\gamma L_g, 4L_f\}$} and $\gamma = \frac{2m\tilde{B}}{\nu}$, and for $\norm{\x_1-\x_\star}^2 \leq B_x$ and $\norm{\y_1 - \y_{j,\star}}^2 \leq B_y$ for all $j$, the output of Algorithm \ref{algo:VRHPS} satisfies
	\begin{align}
		\E{\norm{\x_T - \x_\star}^2} &\leq \col{\tfrac{4\tilde{L}^2}{\mu^2T^2}}(B_x+\tfrac{2(mB_y+2nL_f^2B_x)}{\tilde{L}^2}) \nonumber\\
		&=\O{\tfrac{m^2+m+n}{T^2}}
	\end{align}
	which translates to an SFO complexity of $\O{n+\frac{m+\sqrt{m+n}}{\sqrt{\epsilon}}}$.
\end{theorem}
\begin{IEEEproof}
	Substituting  $\z_t =\x_t - \eta_t\v_t -\eta_t \y_t +\eta_t \y_{j_t,t}$ and $\w_t = \x_{t+1}$ in \eqref{eqlem0} and taking expectation, we obtain
	\begin{align}
		&\En{\norm{\x_{t+1} - \x_\star}^2} \leq \En{\norm{\x_t - \x_\star}^2} \\
		& - (1 - \eta_t \gamma L_g)\En{\norm{\x_{t+1} - \x_t}^2}  \nonumber\\ 
		&  + 2\eta_t \E{\ip{\v_t + \y_t - \y_{j_t,t}}{\x_\star - \x_t} + \phi_{j_t}(\x_\star) - \phi_{j_t}(\x_{t+1})} \nonumber\\
		&  + 2\eta_t \E{\ip{\v_t+ \y_t - \y_{j_t,t}}{\x_t - \x_{t+1}}}\nonumber\\
		&\leqtext{\eqref{hps-proof2}}  \En{\norm{\x_t - \x_\star}^2}  - (1 - \eta_t \gamma L_g)\En{\norm{\x_{t+1} - \x_t}^2}  \nonumber\\ 
		&  + 2\eta_t \E{\ip{\nabla f_{i_t}(\x_t)-\nabla f_{i_t}(\x_\star)}{\x_\star - \x_t}}  \nonumber\\
		&  + 2\eta_t \E{\ip{\v_t + \y_t - \y_{j_t,t} + \nabla \phi_{j_t}(\x_\star)}{\x_t - \x_{t+1}}}\label{vrhps-proof1}
	\end{align}
	where we have also used the unbiased properties of the (sub)gradient estimates to write $\E{\v_t+\y_t-\y_{j_t,t}} = \E{\nabla f_{i_t}(\x_t)}$, so that
	\begin{align}
		&\E{\ip{\v_t + \y_t - \y_{j_t,t} + \nabla \phi_{j_t}(\x_\star)}{\x_\star - \x_t}} \nonumber\\
		&= \E{\ip{\x_\star - \x_t}{ \nabla f_{i_t}(\x_t)  + \nabla \phi(\x_\star)}}  \nonumber\\
		&\eqtext{\eqref{optgrad2}} \E{\ip{\x_\star-\x_t}{ \nabla f_{i_t}(\x_t) - \nabla f(\x_\star)}}\label{svexpeq}
	\end{align}
	\col{The third term on the right of \eqref{vrhps-proof1} can be bounded using Assumption \ref{sc} as in \eqref{hps-t1}. Adding \eqref{vrhps-proof1} with \eqref{eqlem1}, we therefore obtain
    \begin{align}
		&\En{\norm{\x_{t+1} - \x_\star}^2} + 2\eta_t^2\sum_{j=1}^m\En{\norm{\y_{j, t + 1} - \y_{j, \star}}^2} \label{vrhps-proof2}\\
		& \leq(1- \mu \eta_t)\En{\norm{\x_t - \x_\star}^2} + 2\eta_t^2\sum_{j=1}^m \En{\norm{\y_{j,t}-\y_{j,\star}}^2} \nonumber\\
		& - (\tfrac{1}{2} - \eta_t \gamma L_g)\En{\norm{\x_{t+1} - \x_t}^2}  + 2\eta_t^2\En{\norm{\v_t - \nabla f(\x_\star)}^2}\nonumber\\ 
		&  - 2\eta_t \En{D_f(\x_t,\x_\star)}.    \nonumber
	\end{align}
    Here, we note from \eqref{young} that
	\begin{align}
		&\EE\norm{\v_t - \nabla f(\x_\star)}^2 \leq 2\EE\norm{\nabla f_{i_t}(\x_t) - \nabla f_{i_t}(\x_\star)}^2 \nonumber\\
		&+2\EE\norm{\nabla f_{i_t}(\x_\star) - \nabla f_{i_t}(\bx_t) -(\nabla f(\x_\star) - \nabla f(\bx_t))}^2 \nonumber\\
		&\leq 2\EE\norm{\nabla f_{i_t}(\x_t) - \nabla f_{i_t}(\x_\star)}^2 \nonumber\\
		&+ \tfrac{2}{n}\sum_{i=1}^n\EE\norm{\nabla f_i(\bx_t)-\nabla f_i(\x_\star)}^2 	\label{vtdiff}\\
        &\leqtext{\eqref{bound_sm}} 4L_f\EE D_f(\x_t,\x_\star) + \tfrac{4L_f}{n}S_t \label{vtdiff2}
	\end{align}
	where \eqref{vtdiff} follows from the inequality $\EE\norm{\mathsf{X} - \E{\mathsf{X}}}^2 \leq \EE\norm{\mathsf{X}}^2$ and $S_t := \sum_{i=1}^n \E{D_{f_i}(\bx_t, \x_\star)}$. A recursion for $S_t$ can be obtained by observing that $\bx_{t+1}$ is $\x_t$ with probability $1/n$ but remains $\bx_t$ with probability $1-1/n$, so that
    \begin{align}
        S_{t+1} = \EE D_f(\x_t,\x_\star) + \left(1-\tfrac{1}{n}\right)S_t. \label{rec2}
    \end{align}
    Multiplying \eqref{rec2} by $4L_f$ and adding with \eqref{vtdiff2}, we obtain
	\begin{align}
		&4L_fS_{t+1} + \En{\norm{\v_t - \nabla f(\x_\star)}^2} \leq 4L_fS_t + 8 L_f\EE D_f(\x_t,\x_\star). \label{vtrec}
	\end{align}
    Let us define $\Phi_t := 2\sum_{j=1}^m\En{\norm{\y_{j, t} - \y_{j, \star}}^2}
	+ 8L_f S_t$ so that substituting \eqref{vtrec} into \eqref{vrhps-proof2} yields:
	\begin{align}
		\En{\norm{\x_{t + 1} - \x_\star}^2} +& \eta_t^2\Phi_{t+1} \leq (1-\mu \eta_t)\En{\norm{\x_t - \x_\star}^2} + \eta_t^2\Phi_t \nonumber\\
		& - (\tfrac{1}{2} - \eta_t \gamma L_g)\En{\norm{\x_{t+1} - \x_t}^2} \nonumber\\
        &- 2\eta_t\left(1 - 8\eta_tL_f\right)\EE D_f(\x_t,\x_\star). \label{mainrec-vr}
	\end{align}
    The negative terms can be dropped for $\eta_t \leq \frac{1}{\tilde{L}}$ where $\tilde{L} = 2\max\{\gamma L_g, 4L_f\}$. For $\eta_t = \frac{2}{\mu t+2L_f} \leq \tfrac{1}{\tilde{L}}$, we have from \eqref{rect} that 
    }
    \begin{align}
		\E{\norm{\x_T - \x_\star}^2} &\leq \col{\frac{4\tilde{L}^2}{\mu^2T^2}}(B_x+\tfrac{2(mB_y+2nL_f^2B_x)}{\tilde{L}^2}) \nonumber\\
		&\leq \O{\frac{n\kappa_f^2 + m^2\kappa_g^2 + m/\mu^2}{T^2}}\label{vr-rate}
	\end{align}
	where $\O{\cdot}$ hides $B_x$, $B_y$, $\tilde{B}$, $\nu$, and universal constants. To calculate the SFO complexity, we see that there are 3 calls per-iteration on an average and a full gradient evaluation at the initialization. Hence,  \eqref{vr-rate} translates to an SFO complexity of \col{$\O{n+\frac{m+\sqrt{n}}{\sqrt{\epsilon}}}$}.   
    \end{IEEEproof}
We see that the rate in Theorem~\ref{vrhps-thm} strictly improves over the $1/T$ behavior of HPS and the other state-of-the-art algorithms \cite{singh2024stochastic_movingball, singh2024stochastic, singh2025stochastic, xu2020primal, singh2024mini, singh2024unified, necoara2022stochastic}. As expected, variance reduction eliminates the $G_\star$- and $\sigma_\star$-dependent terms from the SFO complexity in Theorem~\ref{hps-thm}. However, in this case, the complexity is still worse than the best-known $\O{\log(1/\epsilon)}$ bound reported in \cite{sanyal2025stochastic} under a different oracle, where all constraint gradients are accessed at every iteration. 

Moreover, accelerated variance-reduction techniques such as Katyusha and its variants, applied only to the smooth component $f_i$ are unlikely to yield further improvements in this constrained setting, since acceleration does not affect either of the main rate-bottlenecks: (a) the factor $(1/2-\eta_t\gamma L_g)$ multiplying $\En{\norm{\x_t - \x_{t+1}}^2}$, and (b) the $\y_{j,t}$-dependent terms required for the telescoping sum. If variance-reduction is not used for the smooth part, the rate becomes $\O{\frac{1}{\epsilon} + \frac{m}{\sqrt{\epsilon}}}$, which is independent of $n$ but still grows linearly with $m$.


In the next section, we switch to the linear regularity assumption \ref{regular} and derive complexity bounds that depend on $r$ rather than explicitly on $m$. Such bounds are particularly useful when the constraints are well-conditioned, so that $r$ stays bounded away from zero and the resulting SFO complexity becomes effectively independent of $m$. 


\section{Nested Hinge-Proximal SGD} \label{sec:n-hps}
In this section, we consider the regime where both $n$ and $m$ are large and seek complexity bounds that are effectively independent of both quantities. To remove the explicit dependence on $n$, we revert to vanilla SGD on the smooth component, without variance reduction. However, as discussed earlier, dependence on $m$ is unavoidable without additional assumptions on the geometry of the constraint region. We therefore rely on the linear-regularity assumption~\ref{regular}, which ensures that the distance of any point to the overall feasible set $\cX$ is controlled by the mean-squared distance to the individual constraint sets $\cX_j$~\cite{patrascu2021stochastic}.
 
A second challenge comes from the penalty parameter $\gamma$ in the exact-penalty reformulation, which typically scales as $\O{m}$ and thus propagates an $m$-dependence into the complexity results. To address these issues, we consider a nested scheme, where the outer loop performs a standard SGD step $\z_t = \x_t - \eta_t \nabla f_{i_t}(\x_t)$ whereas the inner loop computes a constraint-corrected proximal point $\y_t$ for a randomly selected constraint $g_{j_t}$ by solving 
    \begin{align}\label{zsproblem}
        \y_t = \arg\min_{\u \in \Rn^d} &\tfrac{1}{2\eta_t}\norm{\z_t - \u}^2 + h(\u) \nonumber \\ 
        &\text{s.t} \hspace{2mm} g_{j_t}(\u) \leq 0
    \end{align}
    for a randomly selected $j_t$.  To ensure tractability, we utilize the exact penalty reformulation of \eqref{zsproblem}
    \begin{align} \label{e0}
        \y_t = \arg\min_{\u \in \Rn^d} \tfrac{1}{2\eta_t}\norm{\z_t - \u}^2 + h(\u) + \gamma_t[g_{j_t}(\u)]_+
    \end{align}
	which is solved using the hinge-proximal gradient descent updates with step-size $\beta_t\eta_t$:
\begin{align}
    \u_{s+1} &= \prox{\beta_t\eta_t\tphi_{j_t}(\cdot,\u_s)}{(1-\beta_t)\u_s + \beta_t\z_t} \nonumber\\
    &=\arg\min_{\u} \tfrac{1}{2\beta_t\eta_t}\norm{(1-\beta_t)\u_s + \beta_t\z_t-\u}^2 + h(\u) \nonumber\\
    &\qquad + \gamma_t [g_{j_t}(\u_s) + \ip{\nabla g_{j_t}(\u_s)}{\u - \u_s}]_+ \label{nhpsupdate}
\end{align}
As earlier, we define $\phi_{j_t}(\x) = h(\x) + \gamma_t [g_{j_t}(\x)]_+$ and $\tphi_{j_t}(\u,\x) = h(\u) + \gamma_t[g_{j_t}(\x) + \ip{\nabla g_{j_t}(\x)}{\u - \x}]_+$. We can select $\gamma_t$ adaptively to ensure that the problems \eqref{zsproblem} and \eqref{e0} are equivalent. Specifically, applying a standard duality argument (see e.g. \cite[Sec. II-A]{sanyal2025stochastic}) we can bound the dual optimal variable $\lambda_t$ associated with the constraint of \eqref{zsproblem} as 
\col{\begin{align}
	\lambda_t &\leq \frac{\frac{1}{2\eta_t}\norm{\z_t - \tx}^2 + h(\tx) -\frac{1}{2\eta_t}\norm{\z_t - \y_t}^2 - h(\y_t)}{\nu} \\
    &\leq \frac{\norm{\z_t - \tx}^2 + 2\eta_t h(\tx)}{2\eta_t\nu}.
\end{align}}
Since $\z_t$ is known before the inner loop commences, we can simply pick \col{$\gamma_t = \frac{\norm{\z_t - \tx}^2 + 2\eta_t h(\tx)}{2\eta_t\nu}$}, which will ensure that the solution of \eqref{e0} is the same as that of \eqref{zsproblem}. 

Observe that compared to HPS and VR-HPS, the nested scheme above decouples the choice of the penalty parameter $\gamma_t$ from the number of constraints $m$. As a result, $\x_\star$ is no longer characterized by \eqref{optgrad2}; instead, the analysis in this section is based on the optimality condition \eqref{optgrad0}. The full nested hinge-proximal scheme is summarized in Algorithm~\ref{algo:NHPS}. \col{Interestingly, for the special feasibility case $f=h=0$, and for the choice $\beta_t=\tau_t=1$ with sufficiently large $\gamma_t$, Algorithm~\ref{algo:NHPS} reduces to the single-constraint subgradient projection update used in \cite{nedic2019random_feas}.}

\begin{algorithm}
    \caption{Nested HPS}
    \label{algo:NHPS}
    \begin{algorithmic}[1]
        \State \textbf{Parameters:} $\{\eta_t\}$, $\{\gamma_t\}$, $\{\beta_t\}$
        \State \textbf{Input:} initial point $\x_1$        
        \For {$t = 1, \ldots, T-1 $}
        \State Sample $i_t$ and $j_t$ randomly
        \State $\z_t = \x_t - \eta_t \nabla f_{i_t}(\x_t)$
        \State $\u_1 = \x_t$
        \For {$s = 1, \ldots, \tau_t$ } 
        \State $\u_{s+1} = \prox{\beta_t\eta_t\tphi_{j_t}(\cdot,\u_s)}{(1-\beta_t)\u_s + \beta_t\z_t}$
        \EndFor
        \State $\x_{t + 1} = \u_{\tau_t}$
        \EndFor
    \end{algorithmic}
\end{algorithm}

\col{The random minibatch subgradient methods such as in ~\cite{nedic2019random, necoara2022stochastic} also avoid exact projection onto the functional constraint sets. The key distinction is the form of the correction for the sampled constraint $g_{j_t}$. While~\cite{nedic2019random, necoara2022stochastic} use a  subgradient feasibility correction based on the violation $[g_{j_t}(\x_t)]_+$, Algorithm~\ref{algo:NHPS} approximately solves the $g_{j_t}$--corrected proximal problem~\eqref{zsproblem}, with $h$ retained in the inner problem. The inner loop~\eqref{nhpsupdate} implements this correction using only $\prox{h}{\cdot}$ and gradients of $g_{j_t}$.}

The SFO complexity of Algorithm \ref{algo:NHPS} will be characterized in two steps. We first analyze the deterministic inner loop so as to obtain a bound on the squared distance $\norm{\x_{t+1}-\y_t}^2$ in terms of $\norm{\x_t - \y_t}^2$.
\begin{lemma}\label{contract}
	Under Assumption \ref{a1} and for \col{$\gamma_t = \frac{\norm{\z_t - \tx}^2 + 2\eta_t h(\tx)}{2\eta_t\nu}$}, we have that 
	\begin{align} \label{eq:contract}
		\norm{\x_{t + 1}-\y_t}^2 &\leq (\tfrac{1-\beta_t}{1+\beta_t})^{\tau_t}\norm{\x_t - \y_t}^2 
	\end{align}
	for \col{$\beta_t = \tfrac{2\nu}{2\nu + L_g\norm{\z_t - \tx}^2 + 2\eta_t L_g h(\tx)}$}. 
\end{lemma}

\begin{IEEEproof}
	Using the result of Lemma \ref{lem0} for the update in \eqref{nhpsupdate} with appropriate substitutions, we obtain
	\begin{align}
		&\norm{\u_{s+1} - \y_t}^2 \leq \norm{\u_s - \y_t}^2 - (1-\beta_t\eta_t\gamma_t L_g) \norm{\u_{s+1}-\u_s}^2\nonumber\\
		&+2\beta_t\ip{\u_s - \z_t}{\y_t - \u_{s+1}} + 2\beta_t\eta_t(\phi_{j_t}(\y_t) - \phi_{j_t}(\u_{s+1}))\label{nhps-proof1}
	\end{align}
	where we have used the fact that $\u_s - ((1-\beta_t)\u_s + \beta_t\z_t) = \beta_t(\u_s - \z_t)$. Since the objective of \eqref{e0} is $\tfrac{1}{\eta_t}$-strongly convex, we have that
	\begin{align} \label{zs-sc}
		&\norm{\z_t - \y_t}^2 + 2\eta_t(\phi_{j_t}(\y_t) - \phi_{j_t}(\u_{s+1})) \nonumber \\ 
		& + \norm{\u_{s+1}-\y_t}^2 \leq \norm{\z_t - \u_{s+1}}^2
	\end{align}
	which upon multiplying with $\beta_t$ and adding with \eqref{nhps-proof1} yields
	\begin{align}
		&(1+\beta_t)\norm{\u_{s+1} - \y_t}^2 + \beta_t\norm{\z_t - \y_t}^2 \leq \norm{\u_s - \y_t}^2 \nonumber\\
		& - (1-\beta_t\eta_t\gamma_t L_g) \norm{\u_{s+1}-\u_s}^2\nonumber\\
		&+2\beta_t\ip{\u_s - \z_t}{\y_t - \u_{s+1}} +  \beta_t\norm{\z_t - \u_{s+1}}^2.
	\end{align}
	Writing the cross-term as 
	\begin{align}
		2&\ip{\u_s - \z_t}{\y_t - \u_{s+1}} = \norm{\u_{s+1}-\u_s}^2 + \norm{\y_t - \z_t}^2 \nonumber\\
		&-\norm{\z_t - \u_{s+1}}^2 -\norm{\u_s - \y_t}^2
	\end{align}
	we obtain
	\begin{align}
		&(1+\beta_t)\norm{\u_{s+1} - \y_t}^2 \leq (1-\beta_t)\norm{\u_s - \y_t}^2 \nonumber\\
		& - (1-\beta_t - \beta_t\eta_t\gamma_t L_g) \norm{\u_{s+1}-\u_s}^2 
	\end{align}
	where the negative term can be dropped for $\eta_t \leq \tfrac{1-\beta_t}{\beta_t\gamma_tL_g}$, yielding the inner-loop recursion
	\begin{align}
		\norm{\u_{s+1} - \y_t}^2 \leq (\tfrac{1-\beta_t}{1+\beta_t})\norm{\u_s - \y_t}^2 
	\end{align}
	which yields the desired outer-loop recursion. The required condition $\eta_t \leq \frac{1-\beta_t}{\beta_t \gamma_t L_g}$ is guaranteed if we pick \col{$\beta_t = \tfrac{2\nu}{2\nu + L_g\norm{\z_t - \tx}^2 + 2\eta_t L_g h(\tx)}$}. 
	Intuitively, if $\norm{\z_t - \tx}^2$ is too large, $\beta_t$ will be very close to zero and the rate of convergence of the inner loop will be very slow. We will later characterize the expected number of inner loop iterations. 
\end{IEEEproof}

\begin{theorem}\label{nhps-thm}
        Under Assumptions \ref{a1}, \ref{slater}, \ref{sc}, \ref{sig}, \ref{gstarbound}, \ref{regular} and for $\norm{\x_1-\x_\star}^2 \leq B_x$, the average SFO complexity of Algorithm \ref{algo:NHPS} is $\O{\br{\tfrac{1}{\mu^2\epsilon} + \frac{\kappa_f}{\sqrt{\eps}}}\log(\tfrac{1}{\eps})}$.
\end{theorem}

\begin{IEEEproof}
	Using the definition of $\z_t$ and the fact that the objective of \eqref{zsproblem} is $\tfrac{1}{\eta_t}$-strongly convex, we get 
	\begin{align} \label{nthm-proof1}
		&\norm{\y_t - \x_\star}^2 \leq \norm{\x_t - \x_\star}^2 - \norm{\x_t - \y_t}^2 \nonumber\\
		&+ 2\eta_t\ip{\nabla f_{i_t}(\x_t)}{\x_\star - \y_t} - 2\eta_t(\phi_{j_t}(\y_t) - \phi_{j_t}(\x_\star)) 
	\end{align}
	Since $\gamma_t$ is chosen to be sufficiently large, $\y_t$ is also a solution to \eqref{zsproblem} and $[g_{j_t}(\y_t)]_+ = [g_{j_t}(\x_\star)]_+ = 0$. Hence, using the convexity of $h$ and rearranging \eqref{nthm-proof1}, we obtain
	\begin{align}
		&\norm{\y_t - \x_\star}^2 \leq \norm{\x_t - \x_\star}^2 - \norm{\x_t - \y_t}^2 \nonumber\\
		& +2\eta_t\ip{\nabla f_{i_t}(\x_t) - \nabla f_{i_t}(\x_\star)}{\x_\star - \x_t} (=:T_1) \nonumber\\
		&+2\eta_t\ip{\nabla f_{i_t}(\x_\star) - \nabla f(\x_\star)}{\x_\star - \x_t} (=:T_2)\nonumber\\
		&+2\eta_t\ip{\nabla f_{i_t}(\x_t) - \nabla f(\x_\star)}{\x_t - \y_t} (=:T_3)\nonumber\\
		& + 2\eta_t\ip{\nabla f(\x_\star) + \nabla h(\x_\star)}{\x_\star-\y_t}(=:T_4)\label{nthm-proof2}
	\end{align}
	where we have indicated the different summands by $T_1, T_2, T_3$, and $T_4$. We now take expectation and bound the different terms separately. We note that $T_1$ can be bounded using \col{Assumption \ref{sc}} as in \eqref{hps-t1}, $\E{T_2} = 0$, and $\E{T_3}$ can be bounded using \eqref{fvar} as
 	\col{\begin{align}
		\E{T_3} &\leq \tfrac{1}{4}\E{\norm{\x_t-\y_t}^2} + 8\eta_t^2\sigma_\star^2  + 16\eta_t^2L_f \EE D_f(\x_t,\x_\star).\label{nthm-proof3}
	\end{align}}
	For bounding $\E{T_4}$, we introduce $\pi_{\cX}(\x_t)$ and use the optimality condition of \eqref{mainProb} to obtain 
	\begin{align}
		T_4 &= 2\eta_t\ip{\nabla f(\x_\star) + \nabla h(\x_\star)}{\x_\star-\y_t} \nonumber\\
		&\leq 2\eta_t\ip{\nabla f(\x_\star) + \nabla h(\x_\star)}{\pi_{\cX}(\x_t)-\y_t}\nonumber\\
		&= 2\eta_t\ip{\nabla f(\x_\star) + \nabla h(\x_\star)}{\x_t-\y_t} \nonumber\\
		&\quad +2\eta_t\ip{\nabla f(\x_\star) + \nabla h(\x_\star)}{\pi_{\cX}(\x_t)-\x_t}\nonumber\\
		&\hspace{-5mm}\leqtext{\eqref{young},\ref{a1},\ref{gstarbound}} 4\eta_t^2(1+1/r)B_\star +\tfrac{1}{4}\norm{\x_t - \y_t}^2 + \tfrac{r}{4}\norm{\pi_{\cX}(\x_t)-\x_t}^2 \label{nthm-proof4}
	\end{align}
	Also from the regularity assumption, we have that 
	\begin{align}
		r\norm{\pi_{\cX}(\x_t)-\x_t}^2 &\leqtext{\ref{regular}} \EE_{j_t}{\norm{\pi_{\cX_{j_t}}(\x_t)-\x_t}^2} \leqtext{\eqref{zsproblem}} \EE_{j_t}{\norm{\y_t - \x_t}^2} \label{regbound}
	\end{align}
	where \eqref{regbound} uses the definition of $\pi_{\cX_{j_t}}(\x_t)$ and the fact that $\y_t \in \cX_{j_t}$. Taking expectation in \eqref{nthm-proof4} and substituting \eqref{regbound}, we obtain 
	\begin{align}
		\E{T_4} \leq 4\eta_t^2(1+1/r)B_\star + \tfrac{1}{2}\En{\norm{\x_t-\y_t}^2}
	\end{align}
	Hence, taking expectation in \eqref{nthm-proof2} and substituting the bounds for $\E{T_1}, \E{T_2}, \E{T_3}$, and $\E{T_4}$,  we obtain
    \col{\begin{align}
		&\E{\norm{\y_t - \x_\star}^2} \leq \left(1-\mu\eta_t\right)\En{\norm{\x_t - \x_\star}^2} \nonumber\\
		&- \tfrac{1}{4}\E{\norm{\x_t - \y_t}^2} + 4\eta_t^2(1+\tfrac{1}{r})B_\star + 8\eta_t^2\sigma_\star^2 \nonumber\\
		&- 2\eta_t\left(1 - 8\eta_tL_f\right) \EE D_f(\x_t,\x_\star)\label{nthm-proof5}
	\end{align}}
	where we can drop the last term as it is nonpositive for \col{$\eta_t \leq \tfrac{1}{8L_f}$}. Now we use the result of Lemma \ref{contract} into \eqref{nthm-proof5} as
	\begin{align}
		&\En{\norm{\x_{t+1}-\x_\star}^2} \nonumber\\
		&\leq (1+\omega_t)\E{\norm{\x_{t+1}-\y_t}^2} + \left(1+\tfrac{1}{\omega_t}\right)\En{\norm{\y_t-\x_\star}^2} \nonumber\\
		&\leqtext{\ref{eq:contract},\eqref{nthm-proof5}} \left(1+\tfrac{1}{\omega_t}\right)\col{\left(1-\mu\eta_t\right)}\E{\norm{\x_t - \x_\star}^2} \nonumber\\
		&+\E{\left((1+\omega_t)\left(\tfrac{1-\beta_t}{1+\beta_t}\right)^{\tau_t}- \tfrac{1}{4}\left(1+\tfrac{1}{\omega_t}\right)\right)\norm{\x_t - \y_t}^2} \nonumber\\
		&+ 4\eta_t^2\left(1+\tfrac{1}{\omega_t}\right)(1+\tfrac{1}{r})B_\star+ \left(1+\tfrac{1}{\omega_t}\right)8\eta_t^2\sigma_\star^2 \label{nthm-proof6}
	\end{align}
	where the parameter $\omega_t$ can be chosen by setting
	\col{\begin{align}
		\left(1+\tfrac{1}{\omega_t}\right)\left(1-\mu \eta_t\right) &= \left(1-\tfrac{\mu\eta_t}{2} \right) \nonumber\\
		\Rightarrow ~\omega_t &= \tfrac{2}{\mu \eta_t} - 2 \geq 14
	\end{align}}
	for $\eta_t \leq \tfrac{1}{8L_f}$. Hence, if we take $\tau_t$ sufficiently large so as to drop the negative term in \eqref{nthm-proof6} containing $\norm{\x_t - \y_t}^2$, we obtain the recursion
	\begin{align}
		\E{\norm{\x_{t+1}-\x_\star}^2} \leq \col{\left(1-\tfrac{\mu\eta_t}{2}\right)}\E{\norm{\x_t - \x_\star}^2} \nonumber\\
		+ 5\eta_t^2\left((1+\tfrac{1}{r})B_\star+ 2\sigma_\star^2\right)
	\end{align}
	From \eqref{rect}, it follows that for \col{$\eta_t = \frac{4}{\mu t + 32 L_f}$}, we have the bound
	\begin{align}
		\En{\norm{\x_T-\x_\star}^2} &\leq \col{\tfrac{1024L_f^2B_x}{\mu^2T^2} +\tfrac{80}{\mu^2T}} \left((1+\tfrac{1}{r})B_\star+2\sigma_\star^2\right)\nonumber\\
		&=\O{\tfrac{1}{\mu^2 T} + \tfrac{\kappa_f^2}{T^2}} \label{nhps-sfo}
	\end{align}
	which is independent of the number of constraints $m$. Hence, the number of SFO calls to an objective stochastic gradient oracle that returns $\nabla f_{i_t}$ is given by $T = \O{\frac{1}{\mu^2\epsilon} + \frac{\kappa_f}{\sqrt{\eps}}}$. However, since each SFO call returns only a single $\nabla g_{j_t}$ while the inner loop requires $\nabla g_{j_t}(\u_s)$ for every $s = 1, \ldots, \tau_t$, the total number of SFO calls is given by $\sum_{t=1}^{T} \E{\tau_t}$.  
    
    Specifically, we must choose $\tau_t$ such that 
	\col{\begin{align}
		\left(\tfrac{1-\beta_t}{1+\beta_t}\right)^{\tau_t} \leq \tfrac{1}{4\omega_t} &=   \tfrac{\mu \eta_t}{8(1-\mu\eta_t)} \Rightarrow \tau_t \geq \frac{\log(\tfrac{8(1-\mu\eta_t)}{\mu \eta_t})}{\log\left(\tfrac{1+\beta_t}{1-\beta_t}\right)}
	\end{align}}	
	We note the following inequalities:
	\col{\begin{align}
		\log(\tfrac{1+\beta_t}{1-\beta_t}) \geq 2\beta_t &= \tfrac{4\nu}{2\nu + L_g \norm{\z_t - \tx}^2 + 2\eta_t L_gh(\tx)}	 \nonumber\\
		\tfrac{8(1-\mu \eta_t)}{\mu \eta_t}  = 2t + 64 \kappa_f - 8 & \leq 64(t+1)(\kappa_f+1)
	\end{align}
	so that if we choose $\tau_t = \tfrac{1}{2}\log\left[64(t+1)(1+\kappa_f)\right]\left(1+\tfrac{L_g}{2\nu}\norm{\z_t - \tx}^2 + \tfrac{\eta_t L_g}{\nu}h(\tx)\right)$
	then the second summand in \eqref{nthm-proof6} will be negative. Since we can write $\norm{\z_t - \tx}^2 \leq 5\norm{\x_t - \x_\star}^2 + 5\norm{\x_\star - \tx}^2 + 5 \eta_t^2 \norm{\nabla f_{i_t}(\x_t) - \nabla f_{i_t}(\x_\star)}^2 + 5 \eta_t^2 \norm{\nabla f_{i_t}(\x_\star) - \nabla f(\x_\star)}^2 + 5 \eta_t^2 \norm{\nabla f(\x_\star)}^2$, 
	taking expectation and using \eqref{young2} and \eqref{optdist}, we obtain
	\begin{align}
		&\Ex{\norm{\z_t - \tx}^2 + h(\tx)} \leq 5(1+L_f^2)\En{\norm{\x_t - \x_\star}^2} \nonumber\\
		&+ \tfrac{10\tilde{B}}{\mu} + 5\sigma_\star^2 + 5 B_\star + h(\tx) = \O{1}\nonumber
	\end{align}
	where we have used the strong convexity of $f$ and the $\O{\tfrac{1}{t}+\tfrac{\kappa_f^2}{t^2}}$ bound obtained for $\En{\norm{\x_t - \x_\star}^2}$ in \eqref{nhps-sfo}. We also use the fact that $h(\tx)$ is bounded from Assumption \ref{slater}.} Therefore, we have that $\E{\tau_t} \leq \O{\log(t)}$ \col{and the total number of SFO calls is given by 
    \begin{align}
        \sum_{t=1}^{T} \E{\tau_t} \leq \O{\sum_{t=1}^T \log(t)} = \O{T\log(T)}.
    \end{align}
    }   
    Hence, the computational complexity of Algorithm is $\O{(\tfrac{1}{\mu^2\epsilon} + \tfrac{\kappa_f}{\sqrt{\epsilon}})\log\left(\tfrac{1}{\epsilon}\right)}$, with constants independent of $m$.
\end{IEEEproof}

Theorem~\ref{nhps-thm} shows that, under the linear-regularity assumption, Nested HPS attains an SFO complexity of the same order as vanilla SGD for unconstrained strongly convex problems, up to condition-number factors and logarithmic terms, and crucially without any explicit dependence on $m$ or $n$. In contrast to HPS and VR-HPS, whose complexity depends on $m$ through the penalty parameter $\gamma$, the nested scheme keeps $\gamma_t$ bounded and uses the geometry constant $r$ to encapsulate the effect of the constraints. 

\col{The complexity of N-HPS also matches that of \cite{nedic2019random, necoara2022stochastic} up to logarithmic factors. However, our analysis only assumes bounded objective function  gradients at $\x_\star$ (Assumption~\ref{gstarbound}), whereas~\cite[Assumptions~1 and 3]{necoara2022stochastic} and ~\cite[Assumption~1]{nedic2019random} require globally bounded objective function and constraint (sub-)gradients, which can be significantly more restrictive, especially for the strongly convex case considered here.}
     
\section{Numerical Experiments}\label{sec-numerical}
In this section, we demonstrate the effectiveness of the proposed algorithms on the robust regression task introduced in Sec. \ref{sec:intro} on synthetic and real datasets, while generating the perturbations $P_j(\cdot)$ from a corruption simulator. We remark that the parameter $\varepsilon$ plays a critical role as it controls the trade-off between robustness and accuracy. A small $\varepsilon$ tightens the constraint, forcing the solution to be robust under larger sets of perturbations, but at the cost of performance on the original data. Conversely, a larger $\varepsilon$ allows more flexibility but may render the solution less robust. For all simulations, $\varepsilon$ was chosen a priori to ensure the feasible set was nonempty.

We compare the proposed algorithms against representative recent stochastic optimization methods that address similar constrained problems, namely the Primal-Dual Stochastic Gradient (PDSG) algorithm  \cite{xu2020primal},  the Stochastic Subgradient Projection (SSP) method \cite{necoara2022stochastic}, and the L-SVRG algorithm \cite{singh2024unified}. Other works mentioned in Tables \ref{tab:sc-constraints1} and \ref{tab:sc-constraints2} are not used for comparison as they were designed for a different setting and not easily modifiable to handle \eqref{mainProb} \cite{singh2024stochastic_movingball, singh2024stochastic, singh2025stochastic} or are known to have suboptimal SFO complexity  \cite{singh2024mini}. All algorithms are implemented in MATLAB R2024a (M1 processor with 16 GB RAM) and initialized identically with the same random seeds. The parameters for each algorithm were selected individually via a grid search, using performance at the end of $10^6$ SFO calls as the selection criterion. \col{The penalty parameter $\gamma$ is common to all algorithms and is selected by grid search in the range $[10^6, 10^{10}]$. For Algorithm~\ref{algo:NHPS} specifically, the remaining parameters $\beta_t$ and $\tau_t$ are also set to fixed values: $\beta_t = 0.5$ throughout, and $\tau_t$ is capped at a maximum of $30$ inner iterations with early stopping when $\norm{\u_{s+1} - \u_s} < 10^{-3}$. As a result, the inner loop averages $2$--$4$ effective steps per outer iteration in practice.} To improve readability, curves are plotted using a moving average with window size of 1000 consecutive SFO calls, applied uniformly to all methods. 


\subsubsection{Robust regression with synthetic data}
To validate the theoretical convergence properties and assess the computational performance of our proposed algorithms, we conduct comprehensive experiments on controlled synthetic datasets. For simplicity, we consider the squared loss $\ell(x,y) = (x-y)^2$ and generate synthetic datasets with two features and an intercept term, resulting in $\x \in \Rn^3$. We take perturbed training data $P_j(\a_i) = \a_i + \boldsymbol{\delta}_{ij}$ with $\boldsymbol{\delta}_{ij} \sim \cN(\mathbf{0}, \sigma_{\text{train}}^2 \mathbf{I}_3)$. We consider three configurations spanning moderate to large-scale settings, with varying $n$, $m$, and $\varepsilon$ values shown in Table \ref{tab:comparison_methods}. For each configuration, we generate $K = 30$ perturbations per training sample. Finally $\varepsilon$ is chosen separately for each configuration to ensure feasibility.  

\begin{figure*}[htb!]
\centering
\begin{subfigure}[b]{0.33\textwidth}
    \centering
    \includegraphics[width=\textwidth]{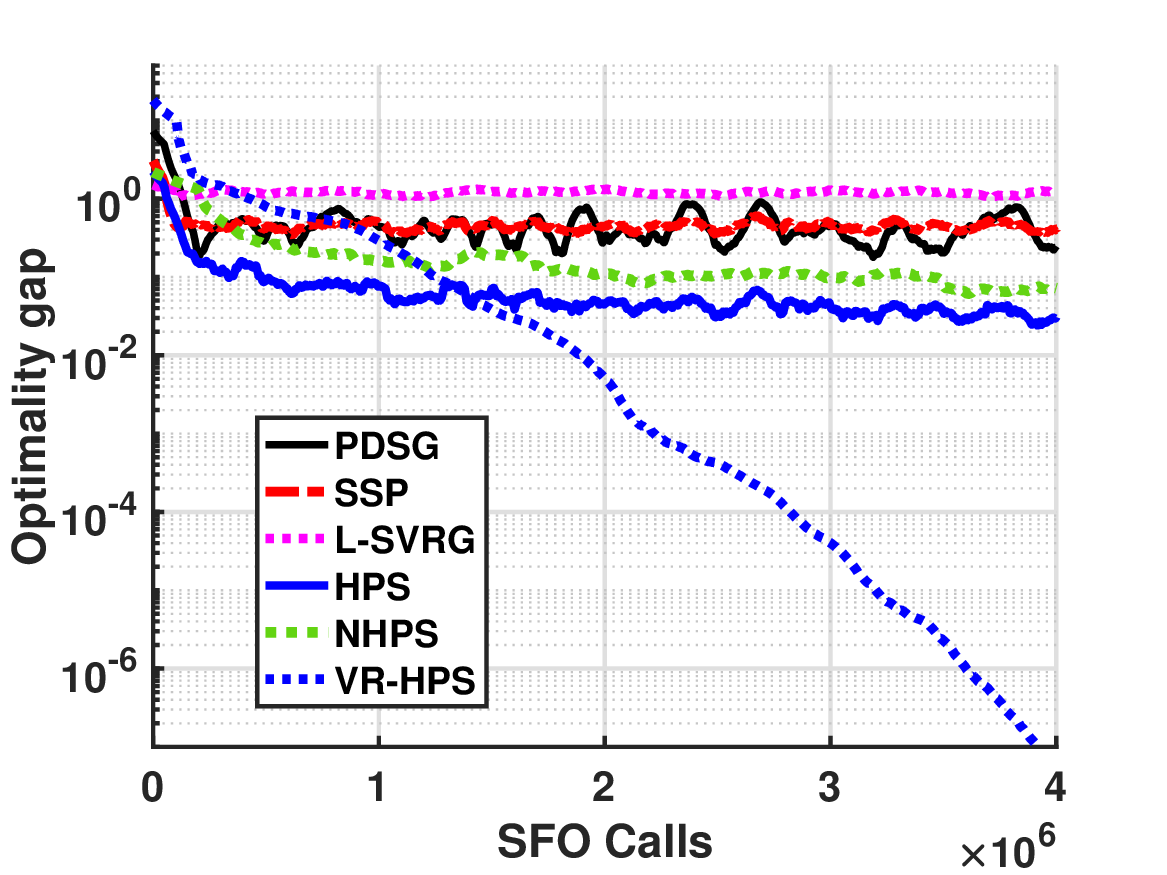}
    \caption{Objective gap $\abs{f(\x_t) - f(\x_\star)}$}
    \label{fig:obj_gap}
\end{subfigure}
\hfill
\begin{subfigure}[b]{0.31\textwidth}
    \centering
    \includegraphics[width=\textwidth]{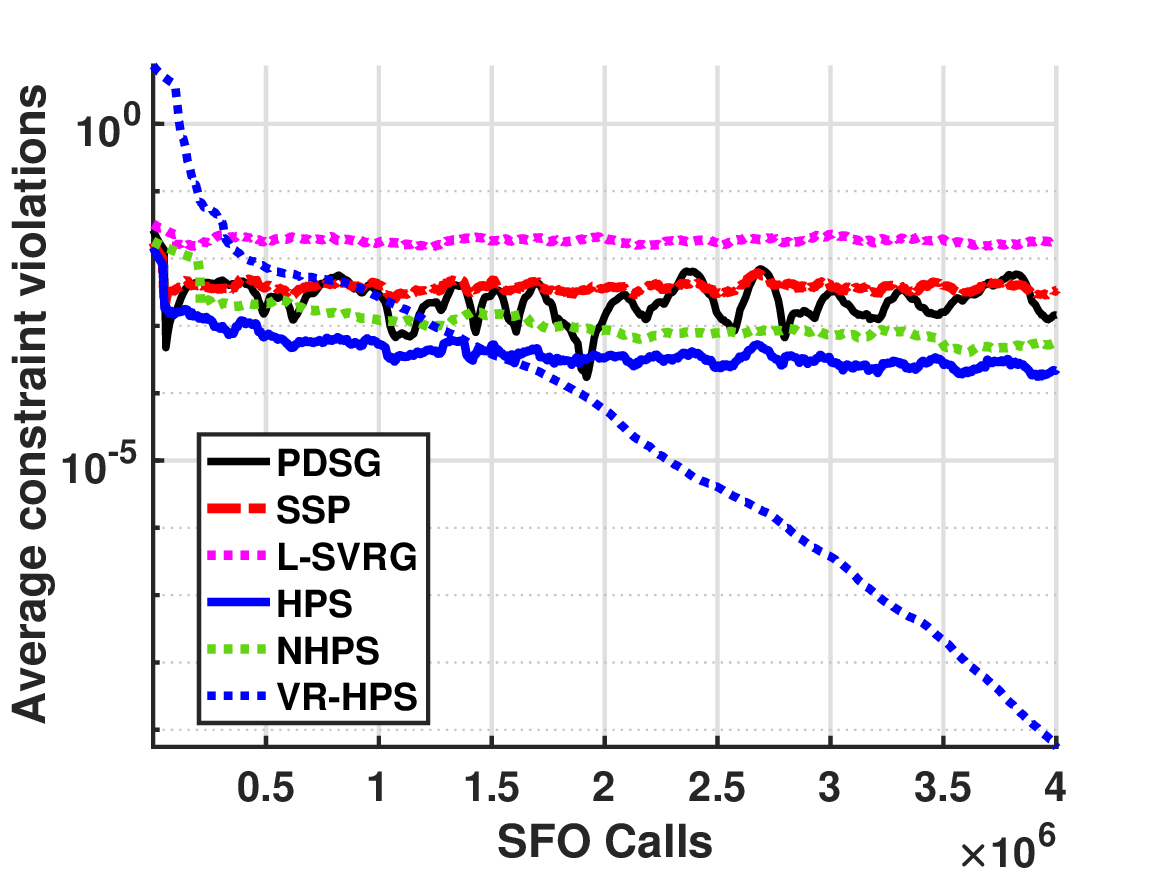}
    \caption{Average constraint violation $\frac{1}{m}\sum_{i=1}^m \bsq{g_i(x_t)}_+$}
    \label{fig:const_vio}
\end{subfigure}
\hfill
\begin{subfigure}[b]{0.33\textwidth}
    \centering
    \includegraphics[width=\textwidth]{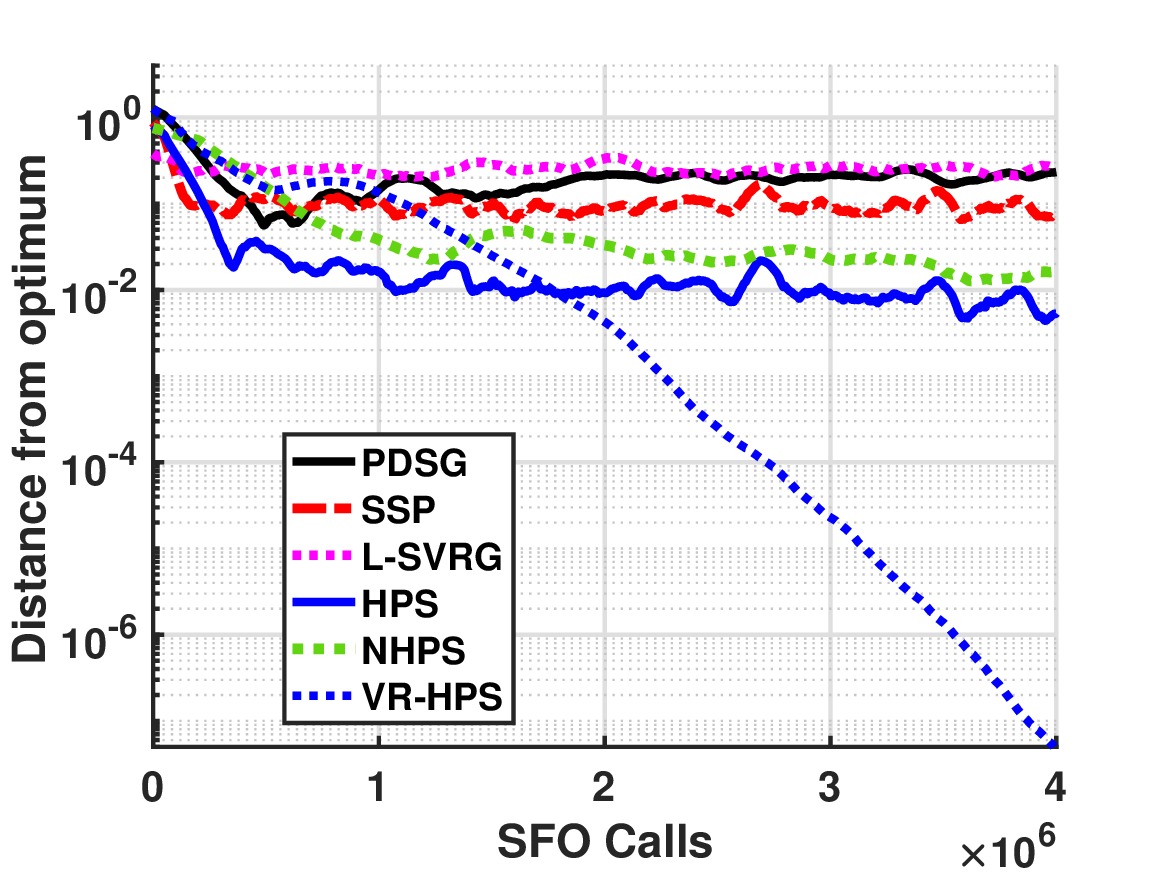}
    \caption{Distance from optimum $\norm{\x_t - \x_\star}$}
    \label{fig:sol_gap}
\end{subfigure}
\caption{\footnotesize{Comparison of proposed stochastic algorithms (HPS, N-HPS) against state-of-the-art baselines (PDSG, SSP), and proposed finite sum algorithm (VR-HPS) against L-SVRG of ~\cite{singh2024unified} on the first configuration with $n=140$, $m=4200$ constraints, and $\varepsilon=13$}}
\label{fig:performance_comparison}
\end{figure*}

Fig. \ref{fig:performance_comparison} plots the evolution of three metrics against the number of SFO calls: (a) objective gap $\abs{f(\x_t) - f(\x^\star)}$, (b) average constraint violation $\frac{1}{m}\sum_{i=1}^m [g_i(\x_t)]_+$, and (c) distance from the optimum $\norm{\x_t - \x^\star}$, for the first configuration with $n = 140$, $m = 4200$ and $\varepsilon = 13$. Among the stochastic algorithms, namely HPS, N-HPS, SSP, and PDSG, we observe that the proposed HPS and N-HPS algorithms perform the best in terms of all three metrics. Among the methods designed to handle finite-sum objectives, namely L-SVRG and VR-HPS, the proposed VR-HPS is clearly and significantly superior. As also predicted by theory, the superiority stems from applying variance reduction to both smooth and non-smooth components of the objective in VR-HPS, as opposed to applying it only to the objective in L-SVRG. 

\begin{table}[htbp]
\caption{\footnotesize{Comparison of VR-HPS against classical methods for robust regression \eqref{eqn:robust_regression}}}
\label{tab:comparison_methods}
\centering
\begin{tabular}{|c|c|c|c|c|c|}
\hline
${n}$ & ${m}$ & $\varepsilon$ & \textbf{Method} & \textbf{RMSE } & \textbf{Wall clock time} (s) \\
\hline
\multirow{7}{*}{140} & \multirow{7}{*}{4200} & \multirow{7}{*}{13} & OLS & 28.102 & 0.007 \\
\cline{4-6}
& & & CVX  & 25.154 & 38.139 \\
\cline{4-6}
& & & DA & 28.433 & 0.005\\
\cline{4-6}
& & & BR & 26.831 &0.004 \\
\cline{4-6}
& & & RF & 28.891 & 3.779\\
\cline{4-6}
& & & EM & 28.433 & 0.011\\
\cline{4-6}
& & & VR-HPS & 25.352 & 2.190\\
\hline
\multirow{7}{*}{350} & \multirow{7}{*}{10500} & \multirow{7}{*}{22} & OLS & 23.84 & 0.008\\
\cline{4-6}
& & & CVX  & 20.132 & 125.088\\
\cline{4-6}
& & & DA & 23.808 & 0.007\\
\cline{4-6}
& & & BR & 22.970 & 0.007 \\
\cline{4-6}
& & & RF & 24.06 & 9.842\\
\cline{4-6}
& & & EM & 23.808 &0.022 \\
\cline{4-6}
& & & VR-HPS & 20.332 & 5.074\\
\hline
\multirow{7}{*}{700} & \multirow{7}{*}{21000} & \multirow{7}{*}{25} & OLS & 25.474 & 0.011\\
\cline{4-6}
& & & CVX  & 21.344 & 874.89 \\
\cline{4-6}
& & & DA & 27.067 & 0.018\\
\cline{4-6}
& & & BR & 26.793 & 0.017\\
\cline{4-6}
& & & RF & 24.549 & 42.011\\
\cline{4-6}
& & & EM &  27.067& 0.074 \\
\cline{4-6}
& & & VR-HPS & 21.484 & 11.955\\
\hline
\end{tabular}
\end{table}

Given the clear theoretical and empirical superiority of VR-HPS, we further evaluate its performance against standard robust regression baselines. To this end, we consider three configurations and utilize 70\% of the generated data (resulting in $n = 140$, $350$, and $700$) to train each algorithm. Subsequently, the test root mean-square error (RMSE) is calculated on the remaining 30\% of the data. Table \ref{tab:comparison_methods} presents comprehensive comparison results between VR-HPS and five baseline approaches (using default parameter settings): Ordinary Least Squares (OLS), Data Augmentation (DA) where the perturbed points are simply added to the data itself, Bayesian Regression (BR), Expectation-Maximization (EM), and Random Forest (RF) regression. Additionally, we also use the CVX optimization package to solve \eqref{eqn:robust_regression}. The results reveal that VR-HPS consistently outperforms all other baselines, coming close to the CVX solution, which also solves \eqref{eqn:robust_regression} optimally. As also evident from the last column of Table \ref{tab:comparison_methods}, for  the largest configuration $(n=700, m = 21000)$, VR-HPS attains an RMSE close to that of CVX while exhibiting almost $73\times$ speedup. The results also suggest that the advantage continues to increase with the problem size. In summary, the proposed robust regression formulation is effective and VR-HPS is an efficient way to solve it. 


\subsubsection{Robust regression with bike sharing data}

\begin{figure*}[!t]
    \centering
    \subfloat[Relative optimality gap $\frac{\abs{f(\x_t) - f(\x_\star)}}{f(\x_\star)}$]{%
        \includegraphics[width=0.4\textwidth]{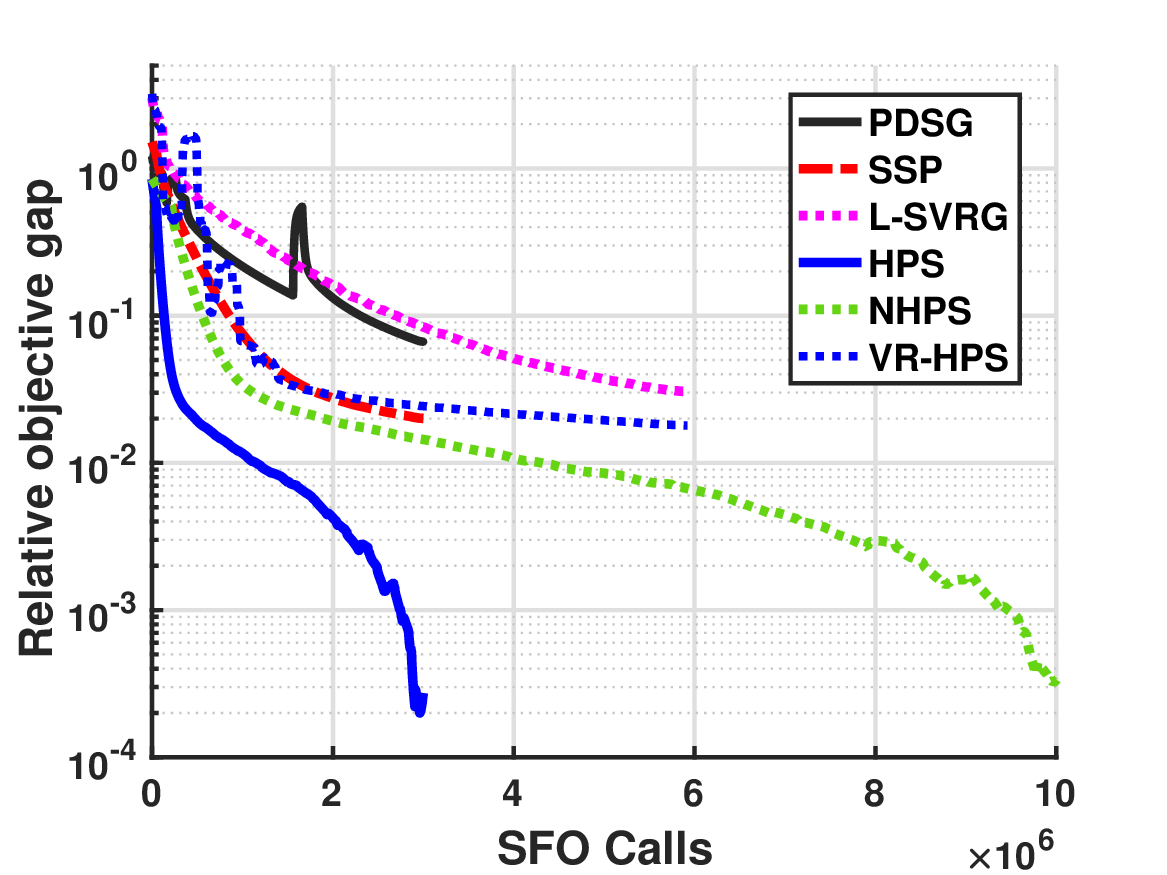}%
        \label{fig:Bike_share_opt_gap}
    }
    \hfil
    \subfloat[Average constraint violation $\frac{1}{m}\sum_{i=1}^m \bsq{g_i(t)}_+$]{%
        \includegraphics[width=0.4\textwidth]{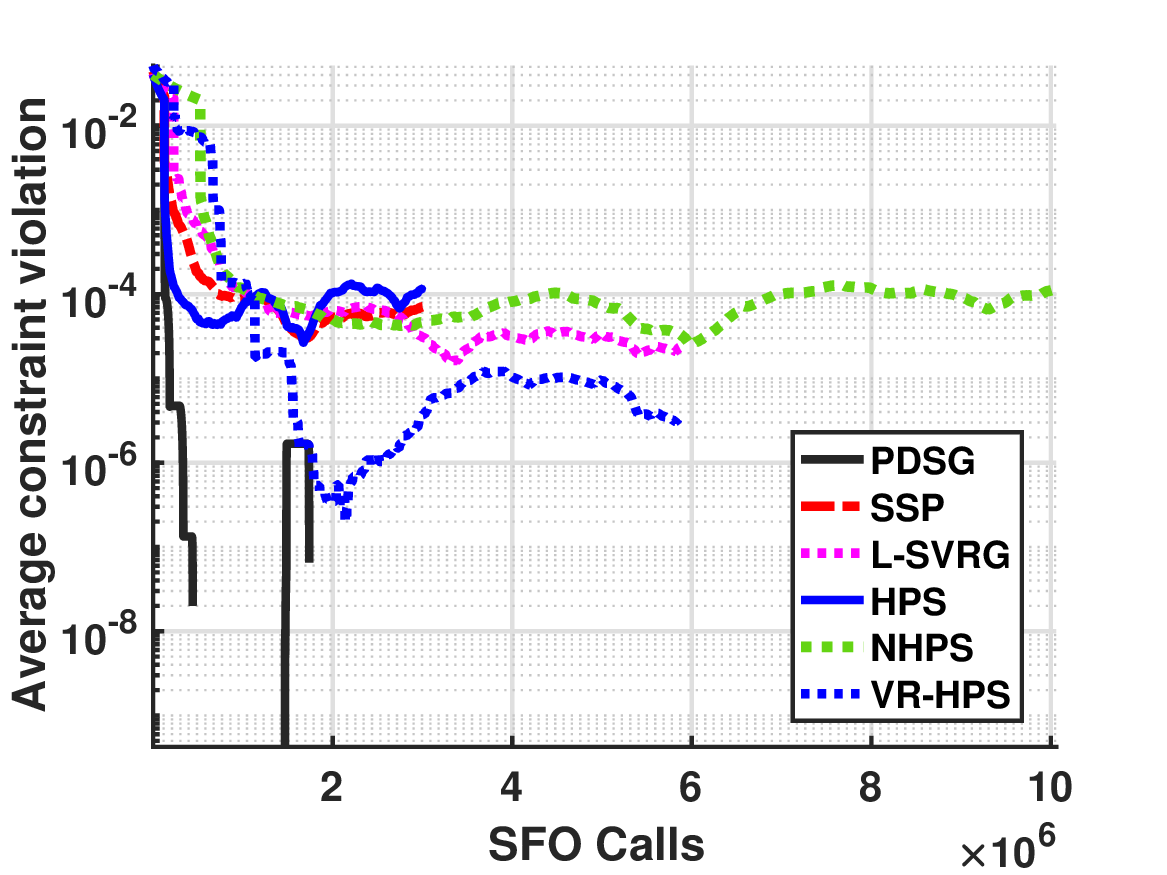}%
        \label{fig:Bike_share_const_vio}
    }
    \caption{\footnotesize{Performance comparison of proposed stochastic algorithms (HPS, N-HPS) against state-of-the-art baselines (PDSG, SSP), and proposed finite sum algorithm (VRHPS) against L-SVRG on Bike-sharing dataset with $m=243300$ constraints. All methods are run for the same iteration budget.}}
    \label{fig:performance_comparison_bikeSharing}
\end{figure*}

The Bike Sharing dataset \cite{FanaeeT2013EventLC}  is a simple yet well-established dataset for benchmarking regression algorithms. The goal here is to predict the number of hourly bike rentals from features such as temporal variables (season, month, hour of day, weekday), weather conditions (temperature, humidity, windspeed), and categorical variables (holiday, working day).  Interestingly, different features exhibit different  reliabilities. Categorical features such as season, year, month, and hour are deterministic and therefore less prone to noise. In contrast, numerical weather-related features (temperature, humidity, and windspeed) are susceptible to sensor measurement errors and forecast inaccuracies. Within the proposed framework in \eqref{eqn:robust_regression}, the perturbed training data $P(\b{a}_i)$ is generated by introducing normally distributed noise to the weather-related features (temperature, humidity, and windspeed) with zero mean and varying standard deviations, while keeping categorical features unperturbed. These perturbations simulate real-world measurement noise as well as weather forecast errors. We use linear model and squared-loss functions for simplicity. 

The experimental setup involves preprocessing the dataset by one-hot encoding the categorical variables and normalizing the numerical features. After one-hot encoding, the feature dimension becomes $d = 51$ including the intercept term. We generate $K = 20$ perturbations per training sample from the $n = 12166$ (70\% of the total samples from the dataset) available training points, resulting in a total of $m = 243300$ constraints. The constraint tolerance is set to $\varepsilon = 550$, chosen through preliminary experiments to balance constraint feasibility with meaningful robustness guarantees. 

Fig. \ref{fig:performance_comparison_bikeSharing} plots the relative optimality gap and the average constraint violation for all algorithms. We observe that in this case, since both $n$ and $m$ are very large, VR-HPS is no longer the best, since it requires $\O{n}$ SFO calls at intermittent steps. Interestingly, for this case, even though $m$ is quite large, the performance of HPS is the best, suggesting that the $\O{m^2}$ bounds obtained in Sec. \ref{sec:proposed_algos} might be loose. Indeed, even if $\gamma$ was increased beyond $10^5$, the performance of HPS did not deteriorate, suggesting that our worst-case $\O{\gamma^2}$ bound, as seen in the proof of Thm. \ref{hps-thm}, may be too conservative. 

The performance of N-HPS is also close but slightly worse than that of HPS, suggesting that while the nested structure eliminates the $m^2$ factor from the SFO complexity of HPS, its empirical performance is worse. Further experiments revealed that the iteration complexity of HPS and N-HPS was almost the same, but N-HPS required, on average, approximately 2.5 inner loops per iteration, which is also evident from Fig. \ref{fig:Bike_share_opt_gap}. In settings where the inner loop is cheap, e.g., if the constraint function gradients can be easily evaluated, the additional cost of N-HPS may be insignificant. 

The test RMSE of HPS for the bike-sharing data was found to be about 3.1\% lower than that of OLS (HPS: $12790$ vs. OLS: $13210$), validating the effectiveness of the robust regression formulation in \eqref{eqn:robust_regression} in a large-scale setting. We remark that the relatively high RMSE is due to the fact that we are using a simple linear model which may not properly fit the real-world data. Nevertheless, this experiment serves to validate the proposed algorithms in a large-scale setting and evaluate their convergence behavior.

\subsubsection{Robust regression with Metro Interstate Traffic Volume data}

\begin{figure*}[!t]
    \centering
    \subfloat[Relative optimality gap $\frac{\abs{f(\x_t) - f(\x_\star)}}{f(\x_\star)}$]{%
        \includegraphics[width=0.4\textwidth]{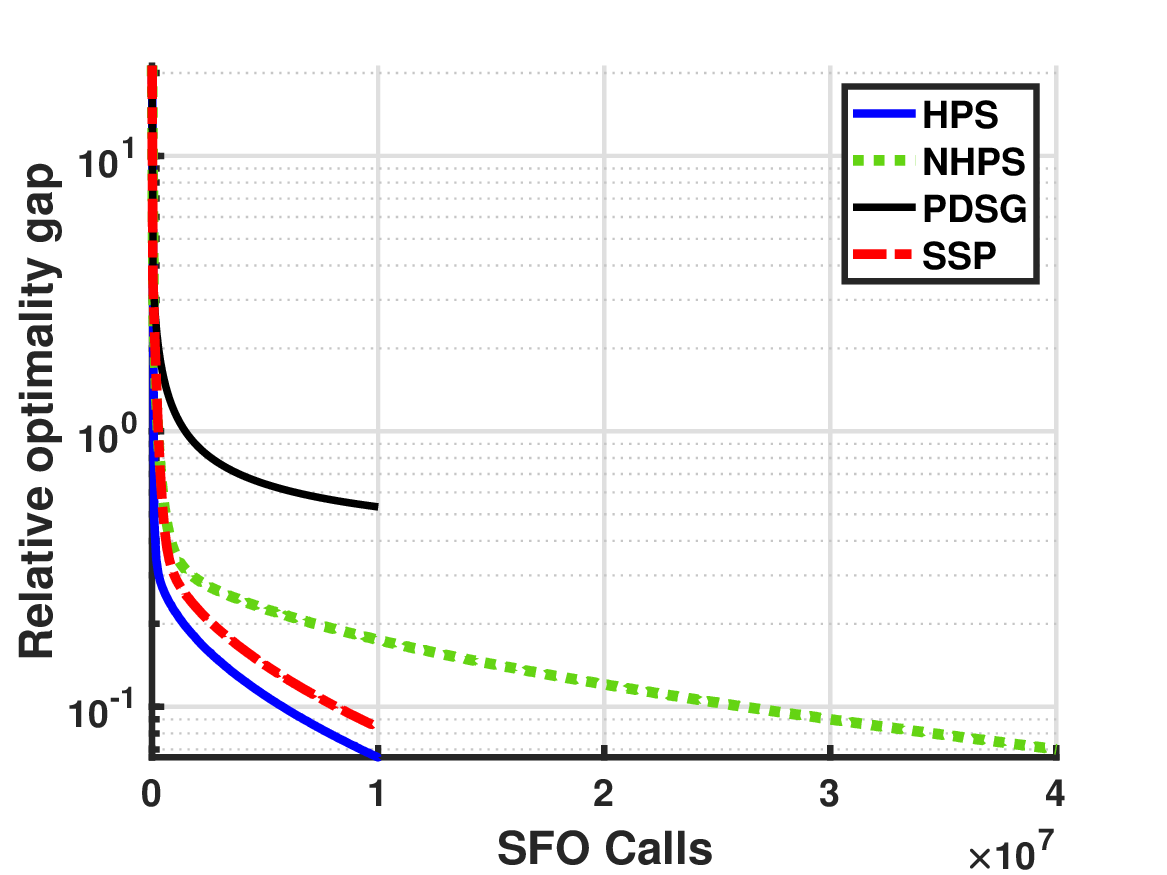}%
        \label{fig:traffic_opt_gap}
    }
    \hfil
    \subfloat[Average constraint violation $\frac{1}{m}\sum_{i=1}^m \bsq{g_i(\x_t)}_+$]{%
        \includegraphics[width=0.4\textwidth]{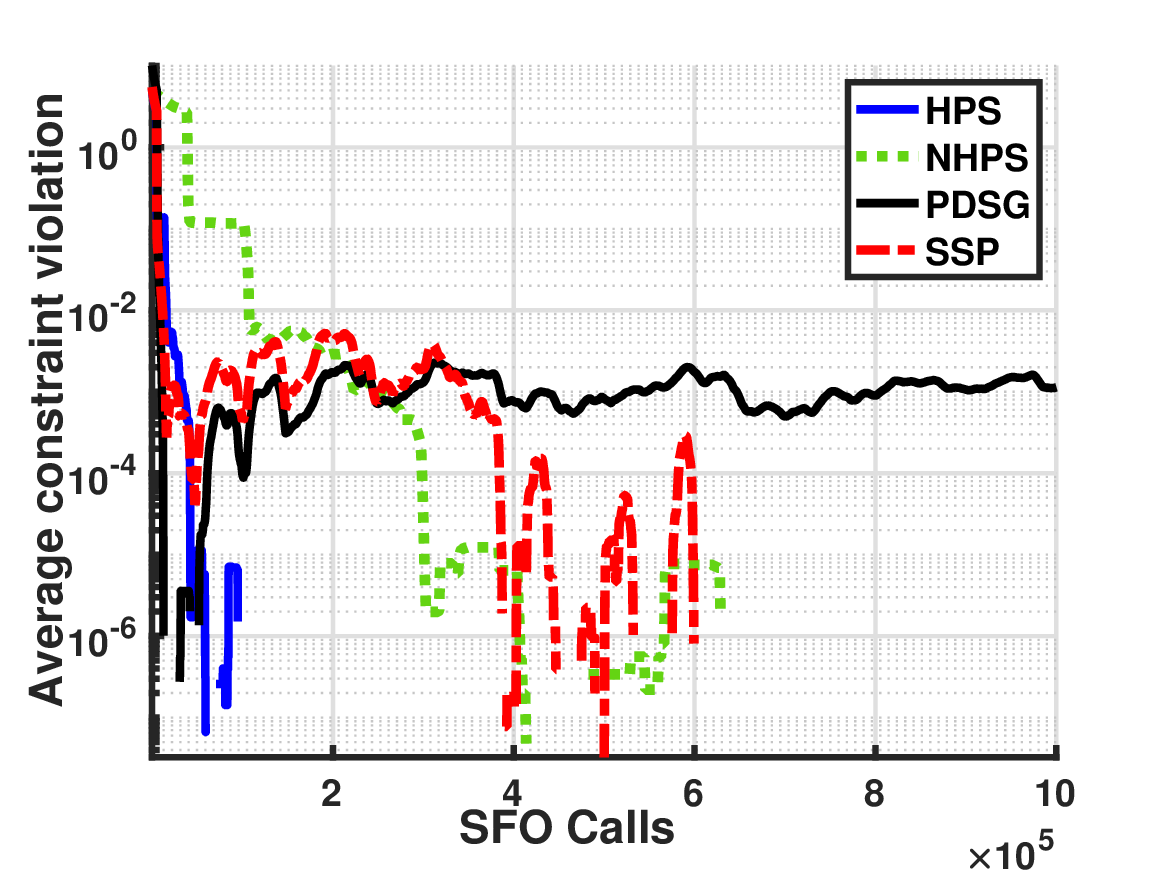}%
        \label{fig:traffic_const_vio}
    }
    \caption{\footnotesize{Performance comparison of HPS, N-HPS, PDSG, and SSP on the Metro Interstate Traffic Volume dataset with $m = 674860$ constraints. All methods are run for the same iteration budget. Fig.~(b) is clipped at $10^6$ SFO calls to show the convergence behavior of all four methods on a common scale; displaying the full range would compress the early phase, where the key differences between methods are most evident.}}
    \label{fig:performance_comparison_traffic}
\end{figure*}

\col{The Metro Interstate Traffic Volume dataset \cite{Hogue2019Metro} contains $48204$ hourly traffic records on the westbound I-94 interstate highway, collected between 2012 and 2018. The goal is to predict hourly traffic volume from features including temporal variables (hour of day, month, day of week), weather conditions (temperature, rainfall, snowfall, cloud cover), holiday indicators, and weather category. The temporal and categorical features are deterministic and hence relatively reliable, whereas continuous weather-related features such as temperature and cloud cover may suffer from sensor errors and forecast inaccuracies. Within the robust regression framework in \eqref{eqn:robust_regression}, the perturbed training data $P(\b{a}_i)$ is generated by adding zero-mean Gaussian noise to the temperature and cloud-cover features with varying standard deviations, while keeping all remaining features unperturbed.}

\col{As part of preprocessing, we one-hot encode the categorical variables: weather category, hour of day, month, and day of week. This results in feature dimension $d=56$, including the intercept term. We use $n = 33743$ training samples (70\% of the total) and generate $K = 20$ perturbations per training sample, yielding a total of $m = 674860$ functional constraints. The constraint tolerance is set to $\varepsilon=3300$, based on preliminary tuning experiments. Due to the large number of constraints, we exclude VR-HPS and L-SVRG since their auxiliary-variable storage scales as $\O{m d}$, which is prohibitive for $m=674860$. We therefore compare only HPS and N-HPS with PDSG and SSP.}

\col{Fig.~\ref{fig:performance_comparison_traffic} shows the relative optimality gap and average constraint violation. In terms of optimality gap, HPS converges fastest with respect to SFO calls, closely followed by SSP. Both these methods reduce the relative optimality gap below $0.1$ within $10^7$ SFO calls. N-HPS is slower in terms of SFO calls, partly because each outer iteration entails multiple inner SFO calls. As compared to the earlier experiments, the convergence of PDSG is significantly slower than that of HPS. Intuitively, for such a large number of constraints and a fixed iteration budget, each dual variable in PDSG is updated infrequently, slowing the accumulation of constraint information in the primal update.}

\col{For constraint violation, HPS, N-HPS, and SSP drive the average violation to near zero within $1\times 10^5$ to $6.5\times 10^5$ SFO calls, after which their curves vanish from the log-scale plot. Together with the bike-sharing experiment, these results show that HPS and N-HPS remain competitive even when the number of constraints is very large.}

\section{Conclusion}\label{sec-conclusion}
We propose a hinge-proximal framework for solving high-dimensional stochastic strongly convex optimization problems with functional constraints. Unlike existing approaches that require globally Lipschitz continuous constraint functions, the proposed framework utilizes exact penalty reformulation that admits smooth constraint functions that need not have bounded gradients. The three instances of the framework, namely hinge-proximal SGD (HPS), variance-reduced HPS, and nested HPS, achieve state-of-the-art or better oracle complexities, while being computationally efficient by using only one constraint gradient per-iteration. We demonstrate the effectiveness of the proposed algorithms on a robust regression problem handling noisy features. 

A critical limitation of our current analysis is that it applies only to strongly convex objectives. It  remains open whether the hinge-proximal framework can be used to remove the bounded-gradient assumption in the general convex (non-strongly convex) and nonconvex regimes. Beyond the three instances proposed here,  the hinge-proximal framework may also be extended to develop distributed variants or projection-free  (Frank–Wolfe–type) algorithms.

	\appendices

    \section{Exact Penalty Reformulation}\label{penaltyequivalent}
Here we prove the equivalence of the original constrained problem \eqref{mainProb} and the reformulated unconstrained problem \eqref{penalty1}.    
We begin by assuming Assumption~\ref{slater} holds for the constrained optimization problem \eqref{mainProb}.  Using the exact penalty method, we rewrite reformulated problem \eqref{penalty1} as 
	\begin{align}
		\x_\star &= \argmin_{\x\in \Rn^d } ~ F(\x) := f(\x) + h(\x) + \frac{\gamma}{m} \sum_{i = 1}^{m}[g_k(\x)]_+\ \nonumber  \\
		& = \argmin_{\x\in \Rn^d, v_i \geq 0} f(\x) + h(\x) + \frac{\gamma}{m} \sum_{i = 1}^{m} v_i \nonumber \\
		& \qquad \quad \text{s.t}~~ g_i(\x) \leq v_i, \hspace{2mm} \forall \nu_i \in \{1,\ldots, m\}  
        \label{penalty1_equivalent}
	\end{align}

For sufficiently large $\gamma$,  the solution of \eqref{penalty1_equivalent} is identical to that of the original problem \eqref{mainProb}. Quantitatively, under Slater's constraint qualification, it suffices to set $\gamma > m\frac{\bt}{\nu}$. To verify this equivalence, we associate dual variables $\mu_k \geq 0$ with the $k$th constraint in \eqref{penalty1_equivalent}, resulting in the Lagrangian as, 
\begin{align}
		\mathcal{L}(\x,v,\mub) &= f(\x) + h(\x) + \frac{\gamma}{m} \sum_{i = 1}^{m} v_i + \sum_{k = 1}^{m} \mu_k (g_k(\x) - v_k) \nonumber \\
		&\hspace{-15mm}= f(\x) + h(\x) + \sum_{k=1}^m \mu_k g_k(\x) + \sum_{i = 1}^{m} v_i(\frac{\gamma}{m} - \mu_i) \nonumber  
	\end{align} 
	where $\mu_i \in \Rn^{m}_+$. Since Slater's condition is satisfied by \eqref{mainProb}, it is also satisfied by \eqref{penalty1_equivalent}. Therefore, the first-order KKT point $(\x_\star,v_\star,\mub_\star)$ satisfies
	\begin{align}
		\hspace{-2mm}(\x_\star,v_\star) =& \argmin_{\x \in \Rn^d, v_i\geq 0} F_h(\x) + \sum_{k=1}^m \mu_{k,\star} g_k(\x) + \sum_{i = 1}^{m} v_i(\frac{\gamma}{m} - \mu_{i, \star})\nonumber
	\end{align}

 where $F_h(\x) = f(\x) + h(\x)$. Hence, for $\frac{\gamma}{m} > \frac{\bt}{\nu} \geq \mu_{i, \star} \hspace{2mm} \forall i \in \{1, \ldots m\}$, it follows that $v_{i, \star} = 0 \hspace{2mm} \forall i \in \{1, \ldots m\}$ and consequently $(\x_\star, \mub_\star)$ is KKT-optimal for \eqref{mainProb}, where $\mub \in \Rn^{m}_+$ collects the dual variables $\{\mu_k\}_{k=1}^m$. 
 Therefore, solving the unconstrained reformulated problem \eqref{penalty1} is equivalent of solving \eqref{mainProb}, provided $\gamma$ is  chosen to be sufficiently large. Hence, we will be using the reformulated problem \eqref{penalty1} in the further convergence analysis of the proposed algorithms. 

\section{Tractability of \eqref{xupdate-sgd}}\label{tractable}
The updates in \eqref{xupdate-sgd} can be written as 
\begin{align}
    &\min_{\u} h(\u) + [\a^\T\u + b]_+ + \frac{1}{2\eta_t}\norm{\z_t - \u}_2^2 \\
    &= \min_{\u,v\geq 0} h(\u) + v + \frac{1}{2\eta_t}\norm{\z_t - \u}_2^2 \label{hlincon}\\ 
    & \text{s. t. } \a^\T\u + b \leq v
\end{align}
where $\a = \gamma\nabla g_j(\x_t)$ and $b = \gamma(g_j(\x_t) -  \x_t^\T \nabla g_j(\x_t))$. Associating dual variable $\lambda$ with the constraint in \eqref{hlincon}, the dual function $\varrho$ is given by 
\begin{align}
    \varrho(\lambda) &= b\lambda + \min_{\u} h(\u) + \lambda\a^\T\u + \tfrac{1}{2\eta_t}\norm{\z_t - \u}_2^2 \nonumber\\
    &\hspace{1cm} + \min_{v\geq 0} v(1-\lambda)  \nonumber\\
    &= b\lambda+\min_{\u} h(\u) + \tfrac{1}{2\eta_t}\norm{\z_t  - \eta_t\lambda\a - \u}_2^2 \nonumber\\
    &\hspace{2cm}- \tfrac{1}{2}\eta_t\lambda^2\a^\T\a + \lambda \a^\T\z_t\nonumber
\end{align}
with the domain $\lambda \in [0,1]$. If $h$ is proximally tractable, we can calculate the optimal $\u$ as $\prox{\eta_t,h}{\z_t-\eta_t\lambda\a}$ and substitute back to find the dual function efficiently. Since the dual problem is a scalar one, we can solve it efficiently and determine $\lambda^\star \in [0,1]$ using a line search algorithm such as the golden section search or bisection, both of which require $O(\log(1/\epsilon))$ iterations. Finally, the solution to \eqref{hlincon} is given by $\prox{\eta_t,h}{\z_t-\eta_t\lambda^\star\a}$. 

In the special case when $h \equiv 0$, we observe that 
\begin{align}
    \varrho(\lambda) = b\lambda - \tfrac{1}{2}\eta_t \lambda^2 \a^\T\a + \lambda \a^\T\z_t \label{closedform}
\end{align}
which is maximized at $\lambda^\star = \min\{\max\{\frac{b+\a^\T\z_t}{\eta_t(\a^\T\a)},0\},1\}$, and the update becomes $\x_{t + 1} = \z_t - \eta_t\lambda^\star \a$. 
    
	\section{Basic inequalities}\label{basic}
	We begin with listing commonly used norm inequalities: for vectors $\u$ and $\v$ and any $\varepsilon > 0$, we have that
	\begin{align}
		\ip{\u}{\v} &\leq \norm{\u}\norm{\v} \leq \tfrac{\varepsilon}{2}\norm{\u}^2+\tfrac{1}{2\varepsilon}\norm{\v}^2 \label{young} \\
		\norm{\u + \v}^2 &\leq \left(1+\tfrac{1}{\varepsilon}\right)\norm{\u}^2 + (1+\varepsilon)\norm{\v}^2 \label{young2}
	\end{align} 
    \col{Since each $f_i$ is smooth, we have the following bound on the Bregman divergence:
    \begin{align}
        D_f(\x,\x_\star) &\geq \frac{1}{2L_f}\EE\norm{\nabla f_{i_t}(\x) - \nabla f_{i_t}(\x_\star)}^2 \label{bound_sm}
    \end{align}
    where the expectation is with respect to the random index $i_t$ for a given $\x$.  Combining with Assumption \ref{sig}, we obtain
    \begin{align}
		&\E{\norm{\nabla f_{i_t}(\x_t) - \nabla f(\x_\star)}^2}\nonumber\\
		&\leq 2\EE\norm{\nabla f_{i_t}(\x_t) - \nabla f_{i_t}(\x_\star)}^2 + 2\EE\norm{\nabla f_{i_t}(\x_\star) - \nabla f(\x_\star)}^2 \nonumber\\
		&\leqtext{\eqref{bound_sm},\eqref{sigeq}} 4L_f \EE D_f(\x_t,\x_\star)
        + 2\sigma_\star^2. \label{fvar}
	\end{align}	
    Strong convexity of $f$ in Assumption \ref{sc} also implies the quadratic lower bound
    \begin{align}
        D_f(\x,\y) &\geq  \frac{\mu}{2}\norm{\x - \y}^2. \label{bound_sc}
    \end{align}
    }

	
	Recall that $\phi_j(\x) = h(\x) + \gamma [g_j(\x)]_+$ and $\nabla h(\x) \in \partial h(\x)$.  Then we see that there exists $\alpha \in [0,1]$, $\nabla h(\x) + \gamma \alpha \nabla g_j(\x) \in \partial \phi_j(\x)$. Hence, for any $\nabla \phi_j(\x) \in \partial \phi_j(\x)$, we have that
	\begin{align}
		&\EE\|\nabla \phi_{j_t}(\x_\star) - \nabla \phi(\x_\star)\|^2 \nonumber\\
		&=\gamma^2\EE\|\alpha_{j_t}\nabla g_{j_t}(\x_\star) - \E{\alpha_{j_t}\nabla g_{j_t}(\x_\star)}\|^2 \nonumber\\
		&\leq \gamma^2\E{\norm{\nabla g_{j_t}(\x_\star)}^2} \leq \gamma^2G_\star^2 \label{phistarbound}
	\end{align}
	We have the following lemma that establishes the convergence of a simple recursion. 
	\begin{lemma}\label{reclem}
		Consider non-negative sequences $\{\sX_t,\sY_t\}$, constants $B$, $L$, and step-size $\eta_t = \frac{2}{at+2L} \leq \frac{1}{L}$ for  $t\geq 0$. Then, we have that
		\begin{subequations}\label{rect}
						\begin{align}
					\sX_{t+1} + \eta_t^2 \sY_{t+1}&\leq (1-a\eta_t)\sX_t + \eta_t^2\sY_t + \eta_t^2B \\
					\Rightarrow 	\sX_T &\leq \tfrac{4L^2}{a^2T^2}(\sX_1+\tfrac{\sY_1}{L^2}) + \tfrac{4B}{a^2T}
				\end{align}
		\end{subequations}
	\end{lemma}
	\begin{IEEEproof}
	We can write the given recursion as
	\begin{align}
		\sX_{t+1} + \tfrac{4}{(at+2L)^2}&\sY_{t+1}\leq \tfrac{at+2L-2a}{at+2L}\sX_t + \tfrac{4B+4\sY_t}{(at+2L)^2} \\
		&\leq \tfrac{(at+2L-a)^2}{(at+2L)^2}\sX_t + \tfrac{4B+4\sY_t}{(at+2L)^2}
	\end{align}
	where we have used the fact that $\tfrac{at+2L-2a}{at+2L} \leq \tfrac{(at+2L-a)^2}{(at+2L)^2}$. Multiplying by $(at+2L)^2$, we obtain the telescopic sum
	\begin{align}
		(at+2L)^2\sX_{t+1} + 4\sY_{t+1}&\leq (at+2L-a)^2\sX_t + 4\sY_t + 4B \nonumber\\
		&\leq4L^2\sX_1 + 4\sY_1 + 4Bt
	\end{align}
	Hence, for non-negative $\sY_t$, we obtain the required bound.
	\end{IEEEproof}
	
	\section{Proof of Lemma \ref{lem1}}\label{prooflem1}
Subtracting $\y_{j_t,\star}$ from \eqref{yjtupdate}, taking norm-square on both sides, and rearranging yields
\begin{align} \label{norm-y}
	& \norm{\y_{j_t, t + 1} - \y_{j_t,\star}}^2 -  \norm{\y_{j_t, t} -  \y_{j_t,\star}}^2 = \tfrac{1}{4\eta_t^2}\norm{\x_t - \x_{t + 1}}^2  \nonumber\\
	& + \norm{\v_t + \y_t}^2 - 2\ip{\y_{j_t, t} -  \y_{j_t,\star}}{\v_t + \y_t} \nonumber \\ 
	& + \tfrac{1}{\eta_t}\ip{\y_{j_t, t} -  \y_{j_t,\star} - \y_t - \v_t}{\x_t - \x_{t + 1}}
\end{align}
Multiplying both sides by $2\eta_t^2$ and introducing $\nabla \phi_{j_t}(\x_\star)$  yields
\begin{align}
	& 2\eta_t^2\norm{\y_{j_t, t + 1} - \y_{j_t,\star}}^2  - 2\eta_t^2\norm{\y_{j_t, t} -  \y_{j_t,\star}}^2 \nonumber\\
	&= \tfrac{1}{2}\norm{\x_t - \x_{t + 1}}^2  + 2\eta_t^2\norm{\v_t + \y_t}^2  \nonumber\\
	&- 4\eta_t^2\ip{\y_{j_t, t} -  \y_{j_t,\star}}{\v_t + \y_t} - T_1 \nonumber 
\end{align}
where $T_1 = 2\eta_t \ip{\v_t +  \y_{j_t,\star} + \y_t - \y_{j_t,t}}{\x_t-\x_{t+1}}$ is the last term in \eqref{eqlem0}. Taking expectation with respect to $j_t$ on both sides, we obtain
\begin{align}
	&\Ej{2\eta_t^2\norm{\y_{j_t, t + 1} - \y_{j_t,\star}}^2  - 2\eta_t^2\norm{\y_{j_t, t} -  \y_{j_t,\star}}^2} \nonumber\\
	&\leq \tfrac{1}{2}\Ej{\norm{\x_t - \x_{t + 1}}^2}  + 2\eta_t^2\norm{\v_t + \y_t}^2  \nonumber\\
	&- 4\eta_t^2\ip{\y_t -  \y_\star}{\v_t + \y_t}- \Ej{T_1} \nonumber \\ 
	&= \tfrac{1}{2}\Ej{\norm{\x_t - \x_{t + 1}}^2}  + 2\eta_t^2\norm{\v_t - \y_\star}^2 - \Ej{T_1} \nonumber\\
	& -2\eta_t^2\norm{\y_t - \y_\star}^2 \label{t1simp}
\end{align}
where we have used the fact that $\x_t$ and $\y_t$ are independent of $j_t$ so that $\Ej{\y_{j_t,t}-\y_{j_t,\star}} = \y_t-\y_\star$. For the term on the left, we proceed as in \cite[Lemma 6]{mishchenko2019stochastic} and write
\begin{align}
	&\Ej{\norm{\y_{j_t, t + 1} - \y_{j_t, \star}}^2} \label{yjt1term}\\
	&= \sum_{j=1}^m\Ej{\norm{\y_{j, t + 1} - \y_{j, \star}}^2} - \Ej{\sum_{j\neq j_t} \norm{\y_{j,t+1}-\y_{j,\star}}^2} \nonumber\\
	&= \sum_{j=1}^m\Ej{\norm{\y_{j, t + 1} - \y_{j, \star}}^2} - \br{1-\tfrac{1}{m}}\sum_{j=1}^m \norm{\y_{j,t}-\y_{j,\star}}^2\nonumber
\end{align}
where the last equality follows from the fact that $\y_{j,t+1}$ for $j\neq j_t$ are not updated at iteration $t$. Hence, we obtain
\begin{align}
	&\Ej{\norm{\y_{j_t, t + 1} - \y_{j_t, \star}}^2 - \norm{\y_{j_t, t} - \y_{j_t, \star}}^2} \nonumber\\
	&=\sum_{j=1}^m\norm{\y_{j, t + 1} - \y_{j, \star}}^2 - \sum_{j=1}^m \norm{\y_{j,t}-\y_{j,\star}}^2
\end{align}
Substituting into \eqref{t1simp}, manipulating, and taking full expectation on both sides, we obtain the required inequality.

\section{Performance of HPS for non-smooth $f$}\label{appendix: HPS_nonsmooth_f}
\col{For the non-smooth case, we first write \eqref{hps-t1} as
\begin{align}
        &\E{\ip{\nabla f_{i_t}(\x_\star) - \nabla f_{i_t}(\x_t)}{\x_\star - \x_t}} \geq \mu \EE\norm{\x_t -\x_\star}^2\label{hps-t1ns}.
\end{align}
where $\nabla f_{i_t}(\x)$ now denotes the subgradient of $f_{i_t}$. Next, assuming that $\EE\norm{\nabla f_{i_t}(\x)}^2 \leq G_f^2$ over the problem domain, we can write \eqref{hps-t2} as
\begin{align}
     & \EE\ip{\nabla f_{i_t}(\x_t) + \nabla \phi_{j_t}(\x_\star)}{\x_t - \x_{t+1}} \nonumber\\
     &\leqtext{\eqref{young}} \tfrac{1}{4\eta_t}\EE\norm{\x_t - \x_{t+1}}^2  + \eta_t\EE\norm{\nabla f_{i_t}(\x_t) + \nabla \phi_{j_t}(\x_\star)}^2 \nonumber\\
    &\leq \tfrac{1}{4\eta_t}\EE\norm{\x_t - \x_{t+1}}^2 + 2\eta_t(4G_f^2 + \gamma^2G_\star^2).\label{hps-t2ns}
\end{align}
which upon combining with \eqref{hps-t1ns} and substituting into \eqref{hps-proof3} yields
\begin{align}
		&\E{\norm{\x_{t+1} - \x_\star}^2} \leq \left(1-2\mu \eta_t \right)\En{\norm{\x_t - \x_\star}^2} \nonumber\\
        &+ 4\eta_t^2(\gamma^2G_\star^2 + 4G_f^2) -(\tfrac{1}{2} - \eta_t \gamma L_g)\En{\norm{\x_{t+1} - \x_t}^2}. \label{hps-proof4ns}
	\end{align}
    Proceeding as in proof of Theorem \ref{hps-thm}, we obtain a similar bound
    \begin{align}
        \EE\norm{\x_T - \x_\star}^2 = \O{\frac{G_f^2+m^2G_\star^2}{\mu^2T} +\frac{m^2\kappa_g^2}{T^2}}.
    \end{align}}

\section{N-HPS for Feasibility Problems}
In this section, we consider solving \eqref{mainProb} for $f = h = 0$. For the specific case of $\tau_t = \beta_t = 1$, the N-HPS algorithm entails carrying out the updates
\begin{align}
    \x_{t+1} = \prox{\eta_t \gamma_t \tg_{j_t}(\cdot,\x_t)}{\x_t}
\end{align}
for $t = 1, \ldots, T$ with $j_t$ being a random index. From \eqref{closedform}, the updates can be written as
\begin{align}
    \x_{t+1} = \x_t - \eta_t \gamma_t\lambda^\star\nabla g_{j_t}(\x_t)
\end{align}
where $\lambda^\star = \pi_{[0,1]}(\frac{[g_{j_t}(\x_t)]_+}{\gamma_t\eta_t\norm{\nabla g_{j_t}(\x_t)}^2})$. Hence, if $\gamma_t > \frac{[g_{j_t}(\x_t)]_+}{\eta_t\norm{\nabla g_{j_t}(\x_t)}^2}$, we obtain the update
\begin{align}
    \x_{t+1} = \x_t - \frac{[\tg_{j_t}(\x_t)]_+}{\norm{\nabla g_{j_t}(\x_t)}^2}\nabla g_{j_t}(\x_t)
\end{align}
\col{which is the subgradient projection onto the linearized halfspace
\[
\{\u:\tg_{j_t}(\u,\x_t)\le 0\}.
\]
Thus, for feasibility problems, N-HPS with $\beta_t=\tau_t=1$ and sufficiently large $\gamma_t$ reduces to the single-sample feasibility update used in~\cite{nedic2019random_feas,nedic2019random}. }

\section{Dependence of SFO complexities on $m$ and $n$} \label{appendix: SFO}
\col{In this section, we analyze the SFO complexities derived in various related works and
explicate their dependence on $m$ and $n$.}

\col{\subsection{Stochastic objectives} For the stochastic objective, the bound provided in \cite[Thm. 3.14]{xu2020primal} takes the form:
\begin{align}
    \EE\norm{\x_{K+1}-\x_\star}^2 \leq \frac{2\alpha}{K+1}\left(\phi_3(\x_\star) + \frac{\log(K+1)}{\rho}\norm{\z_\star}^2\right)
\end{align}
where $\phi_3(\x_\star)$ in \cite[(3.40)]{xu2020primal} contains $C_3$ defined in \cite[(3.37)]{xu2020primal}. Combining these terms, we see that the bound depends in a complicated way on the number of constraints $m$. Careful examination reveals that the dominant term in the bound (for large $m$) is $\norm{\z_\star}^2$ where $\z_\star$ is the dual optimum. Indeed, the augmented Lagrangian
in \cite[(1.2)]{xu2020primal} uses the normalized constraint term $\frac{1}{m}\sum_{j=1}^m \psi_\beta(f_j(\x),z_j)$, and the same normalization appears in \cite[(1.24)]{xu2020primal} and in the KKT condition \cite[(2.1a)]{xu2020primal}. Hence the $\norm{\lam_\star}_1 \leq \bt/\nu = \O{1}$ bound derived from Assumption \ref{slater} in Sec. \ref{sec:assumption} translates to $\norm{\z_\star} \leq \O{m}$. Consequently, the SFO complexity of \cite{xu2020primal} for solving \eqref{mainProb} and under Slater's assumption, is $\Ot{m^2/\epsilon}$. }

\col{\subsection{Deterministic objectives}
The works in \cite{singh2025stochastic, nedic2019random, singh2024stochastic,singh2024stochastic_movingball} are developed for deterministic objectives. When applied to the finite-sum version of \eqref{mainProb}, one evaluation of the deterministic gradient $\nabla f(\x)$ requires $n$ component-gradient evaluations. Therefore, the objective-gradient part of each iteration incurs an additional factor of $n$. In what follows, we convert the iteration complexities of  deterministic-objective methods into SFO complexities using this convention.}

\col{For SGDPA \cite{singh2025stochastic}, we only consider the case where boundedness of
the dual iterates is proved using the perturbation parameter. The perturbed augmented
Lagrangian in \cite[(7)]{singh2025stochastic} uses the normalized term
$\frac{1}{m}\sum_{j=1}^m\psi^j_{\rho,\tau}(h_j(\x),\lambda_j)$, and the KKT condition
\cite[(4)]{singh2025stochastic} contains
$\frac{1-\tau}{m}\sum_{j=1}^m(\lambda_\star)_j\nabla h_j(\x_\star)$. Hence the multiplier
$\lambda_\star$ in \cite{singh2025stochastic} corresponds, up to the fixed factor
$(1-\tau)^{-1}$, to $m\lam_\star$ under the present notation. As before, Slater's
condition gives $\norm{\lambda_\star}=\O{m}$ and therefore
$\norm{\lambda_\star}^2=\O{m^2}$. In the strongly convex case, when the dual iterates
are proved to be bounded, \cite[Thm.~5]{singh2025stochastic} has the bound
$\frac{\norm{\lambda_\star}^2}{m\sqrt K}$ which scales as $\O{m/\sqrt K}$ under Slater's condition. Thus, SGDPA requires
$\Ot{m^2/\epsilon^2}$ iterations and equivalently, an SFO complexity of
$\Ot{nm^2/\epsilon^2}$.}

\col{The works in  \cite{nedic2019random, singh2024stochastic,singh2024stochastic_movingball} use a linear regularity condition that hides away the $m$-dependence. Their strongly-convex rates are all of the order of $1/t$ (keeping the leading terms only), which translates to an SFO complexity of $\O{\frac{n}{\epsilon}}$ for all three. }

\col{\subsection{Finite-sum objectives}
For the mini-batch method of \cite{singh2024mini}, the finite-sum objective and the constraints are both reformulated stochastically; see \cite[(6)--(8)]{singh2024mini}. For partition or nice sampling with objective mini-batch size $\tau_1$ and constraint mini-batch size $\tau_2$, \cite[Thm.~3.5]{singh2024mini} gives $B^2=\frac{n}{\tau_1}\bar B^2$, $L=\frac{n}{\tau_1}\bar L$, and $c=\frac{\bar c m}{\tau_2}$ where we have translated the notation $N$ in \cite{singh2024mini} to $n$ here. The strongly convex rate in \cite[Thm.~4.6]{singh2024mini} is
\[
\EE\norm{\hat \x_K-\x_\star}^2
\leq
\O{\frac{B^2}{\mu^2 C_{\beta,c,B_h}K}},
\]
and for the standard choice $\beta=1$, the constant satisfies
$C_{\beta,c,B_h}=\Theta(1/c)=\Theta(\tau_2/m)$. Hence the iteration complexity is $\O{\frac{nm}{\tau_1\tau_2\epsilon}}$. Since each iteration uses $\tau_1$ objective samples and $\tau_2$ constraint samples in steps \cite[(18)--(20)]{singh2024mini}, the SFO complexity is $\O{(\tau_1+\tau_2)\frac{nm}{\tau_1\tau_2\epsilon}}$. Optimizing this upper bound over $1\leq \tau_1\leq n$ and $1\leq \tau_2\leq m$ yields $\tau_1=\Theta(n)$ and $\tau_2=\Theta(m)$ and the best possible bound $\O{\frac{n+m}{\epsilon}}$. }

\col{For the unified stochastic gradient projection method \cite{singh2024unified}, the
variance-reduced choices are obtained by taking the estimator in Algorithm~3
(SAGA) or Algorithm~4 (L-SVRG). Substituting the corresponding parameter choices into \cite[(19)--(20)]{singh2024unified} gives
$\EE\norm{\hat \x_K-\x_\star}^2=\Ot{1/K}$ where the constants depend only on the regularity parameter in \cite[Assumption~2.5]{singh2024unified}. Hence the
iteration complexity is $\Ot{1/\epsilon}$. Since SAGA requires an initial table/full-gradient
initialization and L-SVRG requires periodic full-gradient refreshes, the corresponding SFO
complexity is $\Ot{n+\frac{1}{\epsilon}}$. }

	\footnotesize
	
	\bibliographystyle{IEEEtran} 
	\bibliography {IEEEabrv,references}


\end{document}